\documentclass[number,citesort,seceqn,dvips]{arxbj}
\usepackage{upgreek}
\usepackage{breakurl}

\aid{0}
\volume{17}
\issue{4}
\pubyear{2011}
\firstpage{1159}
\lastpage{1194}
\doi{10.3150/10-BEJ316}

\makeatletter

\newtheorem{theorem}{Theorem}

\newtheorem{lemma}{Lemma}

\newremark{remark}{Remark}

\makeatother

\begin{document}
\begin{frontmatter}

\title{Multipower variation for Brownian semistationary processes}
\runtitle{Asymptotic results}

\begin{aug}
\author[a]{\fnms{Ole E.} \snm{Barndorff-Nielsen}\thanksref{a}\ead[label=e1]{oebn@imf.au.dk}},
\author[b]{\fnms{Jos\'{e} Manuel} \snm{ Corcuera}\corref{}\thanksref{b}\ead[label=e2]{jmcorcuera@ub.edu}}
\and\\
\author[c]{\fnms{Mark}~\snm{Podolskij}\thanksref{c}\ead[label=e3]{mark.podolskij@math.ethz.ch}}
\runauthor{O.E. Barndorff-Nielsen, J.M. Corcuera and M. Podolskij}
\address[a]{Department of Mathematical Sciences, University of Aarhus,
Ny Munkegade,
DK--8000 Aarhus~C, Denmark. \printead{e1}}
\address[b]{Universitat de Barcelona, Gran Via de les Corts Catalanes
585, 08007
Barcelona, Spain.\\ \printead{e2}}
\address[c]{Department of Mathematics, ETH Z\"{u}rich, HG G32.2, 8092
Z\"{u}rich,
Switzerland.\\ \printead{e3}}
\end{aug}

\received{\smonth{6} \syear{2009}}
\revised{\smonth{6} \syear{2010}}

%
\begin{abstract}
In this paper we study the asymptotic behaviour of power and multipower
variations of processes~$Y$:%
\begin{eqnarray*}
Y_{t}=\int_{-\infty}^{t}g(t-s)\sigma_{s}W(\mathrm{d}s)+Z_{t},
\end{eqnarray*}
where $g\dvtx (0,\infty)\rightarrow\mathbb{R}$ is deterministic, $\sigma
>0$ is a
random process, $W$ is the stochastic Wiener measure and $Z$\ is a
stochastic process in the nature of a drift term. Processes of this type
serve, in particular, to model data of velocity increments of a fluid
in a
turbulence regime with spot intermittency $\sigma$. The purpose of
this paper is to determine the probabilistic limit behaviour of the
(multi)power variations of $Y$ as a basis for studying properties of the
intermittency process $\sigma$. Notably the processes $Y$ are in general
not of the semimartingale kind and the established theory of multipower
variation for semimartingales does not suffice for deriving the limit
properties. As a key tool for the results, a general central limit theorem
for triangular Gaussian schemes is formulated and proved. Examples and an
application to the realised variance ratio are given.
\end{abstract}

%

\end{frontmatter}

\section{Introduction}

The motivation for the development of the results reported in this
paper has
been the need to construct tools for studying the probabilistic limit
behaviour of (realised) quadratic variation and other multipower variations
in relation to the class of \emph{Brownian semistationary} ($\mathcal{BSS}$)
processes. This class, which was introduced in \cite{BNSch09},
consists of
the processes $Y=\{ Y_{t}\} _{t\in\mathbb{R}}$ that are defined
by%
%
\begin{equation}\label{BSS}
Y_{t}=\mu+\int_{-\infty}^{t}g(t-s)\sigma_{s}W(\mathrm{d}s)+\int
_{-\infty
}^{t}q(t-s)a_{s}\,\mathrm{d}s,
\end{equation}
where $\mu$ is a constant; $W$ is a Brownian measure on $\mathbb{R}$; $g$
and $q$ are non-negative deterministic functions on $\mathbb{R}$, with
$%
g( t) =q( t) =0$ for $t\leq0$; and $\sigma$
and $a$
are cadlag processes. When $\sigma$ and $a$ are stationary, so
is $Y$. Hence the name Brownian semistationary processes. It is interesting
to note that the fractional Ornstein--Uhlenbeck process is, in fact, also
representable in the form (\ref{BSS}). The same is true of a wide
class of
stable pseudo-moving average processes; see Corollary~4.3 in~\cite{BNBas10}.

The $\mathcal{BSS}$ processes form the natural analogue, for stationarity
related processes, to the class $\mathcal{BSM}$ of Brownian semimartingales
%
\begin{equation}
Y_{t}=\mu+\int_{0}^{t}\sigma_{s}\,\mathrm{d}W_{s}+\int
_{0}^{t}a_{s}\,\mathrm{d}s.
\end{equation}
In the context of stochastic modelling in finance and in turbulence, the
process $\sigma$ embodies the \emph{volatility} or \emph
{intermittency} of
the dynamics, whether the framework is that of $\mathcal{BSM}$ or
$\mathcal{%
BSS}$. For detailed discussion of $\mathcal{BSS}$ and the more general
concept of tempo-spatial ambit processes see \cite{BNSch04,BNSch07a,BNSch07b,BNSch08a,BNSch08b,BNSch09}. Such
processes are, in particular, able to reproduce key stylized features of
turbulent data.

A main difference between $\mathcal{BSM}$ and $\mathcal{BSS}$ is
that, in
general, models of the $\mathcal{BSS}$ form are not semimartingales
(for a
discussion of this, see Section 3 of \cite{BNSch09}). In consequence,
various important techniques developed for semimartingales, such as the
calculation of quadratic variation by It\^{o} calculus and those of
multipower variation, do not apply or suffice in $\mathcal{BSS}$ settings.
The present paper addresses some of the issues that this raises.

The theory of multipower variation was primarily developed as a basis for
inference on $\sigma$ under $\mathcal{BSM}$ models and, more generally,
It\^{o} semimartingales, with particular focus on inference about the
integrated squared volatility $\sigma^{2+}$ given by%
%
\begin{equation} \label{s2+}
\sigma_{t}^{2+}=\int_{0}^{t}\sigma_{s}^{2}\,\mathrm{d}s.
\end{equation}
This quantity is likewise a focal point for the results discussed in the
following.

Section~\ref{sec2} introduces common notation for multipower variation and recalls
some basic properties of such quantities. A law of large numbers and a
central limit theorem for multipower variation of triangular arrays of
Gaussian random variables are derived in Section~\ref{sec3}, and these limit results
are drawn upon in Section~\ref{sec4} to establish probability and central limit
theorems for multipower variation for $\mathcal{BSS}$\ processes, with most
of the proofs postponed to Section~\ref{sec8}. Section~\ref{sec5} presents several
examples and Section~\ref{sec6} discusses an application concerning the limit
behaviour of the realised variation ratio, that is, the ratio of realised
bipower variation to realised quadratic variation. Section~\ref{sec7} concludes and
indicates some possible directions for further related work.

\section{Multipower variation}\label{sec2}

The concept of (realised) multipower variation was originally
introduced in %
\cite{BNS04a} in the context of semimartingales, and the mathematical theory
has been studied in a number of papers \cite{BNGJPS06,Jac08a,BNSW06,KiPo07} while various applications are the main subjects in %
\cite{BNS04b,BNS06,BNS07,Jac08b,Woe06}.
Multipower variation turns out to be useful for analysing properties of
parts of a process that are not directly observable. In this section we
present the definition of realised multipower variation and recall its
asymptotic properties for some classes of processes.

Let us consider a continuous-time process $X$, defined on some filtered
probability space $(\Omega,\mathcal{F},(\mathcal{F}_{t})_{t\geq0},P)$,
that is observed at equidistant time points $t_{i}=i/n$, $i=0,\ldots,[nt]$.
A realised multipower variation of the process $X$ is an object of the type

\begin{equation}\label{multipowerdef}
\sum_{i=1}^{[nt]-k+1}\prod_{j=1}^{k}|\Delta
_{i+j-1}^{n}X|^{p_{j}},\qquad
\Delta_{i}^{n}X=X_{{i/n}}-X_{{(i-1)/n}},\qquad
p_{1},\ldots
,p_{k}\geq0 ,
\end{equation}
for some fixed number $k\geq1$. We now present an overview of the
asymptotic theory for quantities of the form (\ref{multipowerdef}) for
various types of processes $X$.

We start with the $\mathcal{BSM}$ case
%
\begin{equation}\label{contsm}
X_{t}=X_{0}+\int_{0}^{t}a_{s}\,\mathrm{d}s+\int_{0}^{t}\sigma
_{s}\,\mathrm{d}%
W_{s},
\end{equation}
where $W$ is a Brownian motion, $a$ is a locally bounded and predictable
drift process and $\sigma$ is an adapted and cadlag volatility
process. As was established in \cite{BNGJPS06}, the convergence in
probability
%
\begin{equation}\label{smlln}
n^{p_{+}/2-1}\sum_{i=1}^{[nt]-k+1}\prod_{j=1}^{k}|\Delta
_{i+j-1}^{n}X|^{p_{j}}\stackrel{\mathrm{ucp}}{\longrightarrow}\mu
_{p_{1}}\cdots\mu
_{p_{k}}\int_{0}^{t}|\sigma_{s}|^{p_{+}}\,\mathrm{d}s
\end{equation}
holds, where $p_{+}=\sum_{j=1}^{k}p_{j}$ and $\mu_{p}=E[|u|^{p}]$,
$u\sim
N(0,1)$ and we write $Z^{n}\stackrel{\mathrm{ucp}}{\longrightarrow}Z$ when $%
\sup_{t\in\lbrack0,T]}|Z_{t}^{n}-Z_{t}|\stackrel{P}{\longrightarrow
}0$ for
any $T>0$. Under a further condition on the volatility process, one obtains
the associated stable central limit theorem:
\begin{eqnarray}\label{smclt}
&&\sqrt{n}\Biggl(n^{p_{+}/2-1}\sum_{i=1}^{[nt]-k+1}\prod
_{j=1}^{k}|\Delta
_{i+j-1}^{n}X|^{p_{j}}-\mu_{p_{1}}\cdots\mu_{p_{k}}\int
_{0}^{t}|\sigma
_{s}|^{p_{+}}\,\mathrm{d}s\Biggr)\nonumber
\\[-8pt]\\[-8pt]
&&\quad\stackrel{\mathrm{st}}{\longrightarrow}\sqrt{A}%
\int_{0}^{t}|\sigma_{s}|^{p_{+}}\,\mathrm{d}B_{s} ,\nonumber
\end{eqnarray}
where $B$ is another Brownian motion, defined on an extension of the
probability space $(\Omega,\mathcal{F},(\mathcal{F}_{t})_{t\geq
0},P)$ and
independent of $\mathcal{F}$, and the constant $A$ is given by
\begin{eqnarray*}
A=\prod_{l=1}^{k}\mu_{2p_{l}}-(2k-1)\prod_{l=1}^{k}\mu
_{p_{l}}^{2}+2\sum_{m=1}^{k-1}\prod_{l=1}^{m}\mu
_{p_{l}}\prod_{l=k-m+1}^{k}\mu_{p_{l}}\prod_{l=1}^{k-m}\mu_{p_{l}+p_{l+m}}.
\end{eqnarray*}

Recall that the stable convergence of processes is defined as follows. A
sequence of processes $Z^{n}$ converges stably in law towards the
process $Z$
(written $Z^{n}\stackrel{\mathrm{st}}{\longrightarrow}Z$), that is, defined on the
extension of the original probability space $(\Omega,\mathcal
{F},(\mathcal{F%
}_{t})_{t\geq0},P)$, if and only if for any bounded and continuous
real-valued functional $f$ and any $\mathcal{F}$-measurable random
variable $%
V$ it holds that
\begin{eqnarray*}
\lim_{n\rightarrow\infty}E[f(Z^{n})V]=E[f(Z)V],
\end{eqnarray*}
we use the notation $Z^{n}\stackrel{\mathrm{st}}{\longrightarrow}Z$.

A crucial property of the realised multipower variation is its
robustness to
jumps when $\max_{i}(p_{i})<2$ \cite{BNSW06,Jac08c}.
Assume for a
moment that $X$ is a general It\^{o} semimartingale with continuous
part $%
X^{c}$ satisfying (\ref{contsm}). Then, by (\ref{smlln}) and the robustness
property, we obtain the convergence
\begin{eqnarray*}
\frac{\mu_{1}^{-2}\sum_{i=1}^{[nt]-1}|\Delta_{i}^{n}X||\Delta
_{i+1}^{n}X|%
}{\sum_{i=1}^{[nt]}|\Delta_{i}^{n}X|^{2}}\stackrel{P}{\longrightarrow
}\frac{%
[X^{c}]}{[X]} ,
\end{eqnarray*}
where $[X]$ denotes the quadratic variation of the semimartingale $X$ and
the limit is less than or equal to 1. The latter result, together with the
stable convergence in (\ref{smclt}), can be used to construct a formal test
for jumps (see \cite{BNS04a}). On the other hand, we know that if the limit
of the left-hand side is greater than 1 (which is the case for some typical
turbulence data), the process $X$ cannot be an It\^{o} semimartingale.

In another direction, a study \cite{BNCPW08} was made of the asymptotic
behaviour of bipower variation for processes of the type
%
\begin{equation}\label{contgaus}
X_{t}=X_{0}+\int_{0}^{t}\sigma_{s}\,\mathrm{d}G_{s} ,\qquad t\geq0 ,
\end{equation}
where $G$ is a continuous Gaussian process with centered and stationary
increments (the latter integral is defined as a Riemann--Stieltjes integral).
The process defined in (\ref{contgaus}) is, in general, also not a
semimartingale, and the theory in \cite{BNGJPS06} does not apply. In
particular, a different normalisation is required. Define the (normalised)
multipower variation by
\begin{eqnarray*}
V(X,p_{1},\ldots,p_{k})_{t}^{n}=\frac{1}{n\tau_{n}^{p_{+}}}%
\sum_{i=1}^{[nt]-k+1}\prod_{j=1}^{k}|\Delta
_{i+j-1}^{n}X|^{p_{j}} ,\qquad
p_{1},\ldots,p_{k}\geq0 ,
\end{eqnarray*}
where $\tau_{n}>0$ is given by
%
\begin{equation}\label{tau}
\tau_{n}^{2}=\bar{R}[(1/n)],
\end{equation}
with
%
\begin{equation}\label{Rbar}
\bar{R}(t)=E[(G_{s+t}-G_{s})^{2}].
\end{equation}

Under some assumptions on $\bar{R}$ and the volatility process $\sigma
$ it
was shown that
\begin{eqnarray*}
V(X,p_{1},\ldots,p_{k})_{t}^{n}\stackrel{\mathrm{ucp}}{\longrightarrow}\rho
_{p_{1},\ldots,p_{k}}\int_{0}^{t}|\sigma_{s}|^{p_{+}}\,\mathrm{d}s
\end{eqnarray*}
for a certain constant $\rho_{p_{1},\ldots,p_{k}}$ that depends on the
behaviour of $\bar{R}$ near 0. Furthermore, an associated (stable) central
limit theorem, of a form similar to (\ref{smclt}), was derived. Note, however,
that in general there are essential differences between the characters
of $%
\mathcal{BSS}$ processes and processes of type (\ref{contgaus}). In the
latter case, the process $\sigma$ has only a local influence in  the value
of $X$ whereas, in the $\mathcal{BSS}$ case, the process is also
affected by
the past of $\sigma.$

\section{Multipower variation of Gaussian triangular arrays}\label{sec3}

In this section we derive some asymptotic results for functionals of
arrays of stationary Gaussian sequences. We consider a triangular array
$%
(X_{i,n})_{n\geq1,1\leq i\leq\lbrack nt]}$ ($t>0$) of row-wise stationary
Gaussian variables with mean $0$ and variance $1$. Let
%
\begin{equation}\label{corrfunction}
r_{n}(j)=\operatorname{cor}( X_{1,n},X_{1+j,n})  ,\qquad j\geq0,
\end{equation}
be the correlation function of $(X_{i,n})_{1\leq i\leq\lbrack nt]}$. Assume
that the array $(X_{i,n})_{n\geq1,1\leq i\leq\lbrack nt]}$ is
``non-degenerate'', that is, the covariance matrix of $(X_{i,n},\ldots
,X_{i+k,n})$
is invertible for any $k\geq1$ and $n\geq1$ (otherwise the results below
do not hold).

Now, define the multipower variation associated with the sequence $%
(X_{i,n})_{n\geq1,1\leq i\leq\lbrack nt]}$:
%
\begin{equation}\label{multiarray}
V(p_{1},\ldots,p_{k})_{t}^{n}=\frac{1}{n}\sum_{i=1}^{[nt]-k+1}%
\prod_{j=1}^{k}|X_{i+j-1,n}|^{p_{j}} ,\qquad p_{1},\ldots,p_{k}\geq0.
\end{equation}
Our first result is the weak law of large numbers.

\begin{theorem}
\label{th1} Assume that there exists a sequence $\mathrm{r}(j)$ with
%
\begin{equation}\label{pncond}
r_{n}^{2}(j)\leq\mathrm{r}(j) ,\qquad\frac{1}{n}\sum
_{j=1}^{n-1}\mathrm{r}%
(j)\rightarrow0
\end{equation}
as $n\rightarrow\infty$. Then it holds that
%
\begin{equation} \label{pnconv}
V(p_{1},\ldots,p_{k})_{t}^{n}-\rho_{p_{1},\ldots
,p_{k}}^{(n)}t\stackrel{\mathrm{ucp}%
}{\longrightarrow}0 ,
\end{equation}
where
%
\begin{equation} \label{center}
\rho_{p_{1},\ldots,p_{k}}^{(n)}=E[|X_{1,n}|^{p_{1}}\cdots
|X_{k,n}|^{p_{k}}].
\end{equation}
\end{theorem}

\begin{pf}
 See Section \ref{sec8}.
\end{pf}

Before we present the associated central limit theorem, we need to introduce
another Gaussian process.  Suppose that $r_{n}(j)\rightarrow\rho
(j)$, $%
j=1,\ldots,k-1$, for some numbers $\rho(j).$ Let $(Q_{i})_{i\geq1}$
be a
non-degenerate, stationary, centered (discrete-time) Gaussian process with
variance 1 and correlation function
%
\begin{equation}\label{corrQ}
\rho(j)=\operatorname{cor}( Q_{1},Q_{1+j})  ,\qquad j\geq1.
\end{equation}
Define
%
\begin{equation}\label{VQ}
V_{Q}(p_{1},\ldots,p_{k})_{t}^{n}=\frac{1}{n}\sum_{i=1}^{[nt]-k+1}%
\prod_{j=1}^{k}|Q_{i+j-1}|^{p_{j}}
\end{equation}
and let $\rho_{p_{1},\ldots,p_{k}}=E[|Q_{1}|^{p_{1}}\cdot\cdot\cdot
|Q_{k}|^{p_{k}}].$ Then $\rho_{p_{1},\ldots,p_{k}}^{(n)}\rightarrow
\rho
_{p_{1},\ldots,p_{k}}$ and in this case we obtain the uniform
convergence on compacts in probability (or ucp convergence):
\begin{eqnarray*}
V(p_{1},\ldots,p_{k})_{t}^{n}\stackrel{\mathrm{ucp}}{\longrightarrow}\rho
_{p_{1},\ldots, p_{k}}t.
\end{eqnarray*}
Now we have the following central limit theorem for the family $%
(V(p_{1}^{j},\ldots,p_{k}^{j})_{t}^{n})_{1\leq j\leq d}$.

\begin{theorem}
\label{th2} Assume that
%
\begin{equation}\label{cltcond1}
r_{n}(j)\rightarrow\rho(j) ,\qquad j\geq0,
\end{equation}
and that, for any $j,n\geq1$, there exists a sequence $\mathrm{r}(j)$ with
%
\begin{equation}\label{cltcond2}
r_{n}^{2}(j)\leq\mathrm{r}(j) ,\qquad\sum_{j=1}^{\infty}\mathrm{r}
(j)<\infty.
\end{equation}
Then we have
%
\begin{equation} \label{beta2}
\sqrt{n}\bigl(V(p_{1}^{j},\ldots,p_{k}^{j})_{t}^{n}-\rho
_{p_{1}^{j},\ldots
,p_{k}^{j}}^{(n)}t\bigr)_{1\leq j\leq d}\stackrel{\mathcal
{L}}{\rightarrow}%
\beta^{1/2}B_{t} ,
\end{equation}
where $B$ is a $d$-dimensional Brownian, $\beta$ is a $d\times d$%
-dimensional matrix given by
%
\begin{equation}\label{beta}
\beta_{ij}=\lim_{n\rightarrow\infty}n \operatorname{cov}
(V_{Q}(p_{1}^{i},%
\ldots,p_{k}^{i})_{1}^{n},V_{Q}(p_{1}^{j},\ldots
,p_{k}^{j})_{1}^{n})%
 ,\qquad1\leq i,j\leq d ,
\end{equation}
and the weak convergence holds in the space $\mathcal{D}([0,T]^{d})$
equipped with the uniform topology.
\end{theorem}

\begin{pf} See Section~\ref{sec8}.\end{pf}

\begin{remark}
Similar asymptotic results can be obtained for general quantities of the
form
%
\begin{equation}\label{generalmulti}
\frac{1}{n}\sum_{i=1}^{[nt]-k+1}H(X_{i,n},\ldots,X_{i+k-1,n})
\end{equation}
for some function $H\dvtx \mathbb{R}^{k}\rightarrow\mathbb{R}$. Let $m$
denote the
Hermite index of $H$ (notice that the Hermite index of the power function
used in (\ref{multiarray}) is 2). Replace condition (\ref{pncond}) by
\begin{eqnarray*}
|r_{n}^{m}(j)|\leq\mathrm{r}(j) ,\qquad\frac{1}{n}\sum
_{j=1}^{n-1}\mathrm{r%
}(j)\rightarrow0
\end{eqnarray*}
and (\ref{cltcond2}) by
\begin{eqnarray*}
|r_{n}^{m}(j)|\leq\mathrm{r}(j) ,\qquad\sum_{j=1}^{\infty}\mathrm
{r}%
(j)<\infty.
\end{eqnarray*}
Then Theorems \ref{th1} and \ref{th2} hold true for the functional
(\ref%
{generalmulti}) provided that $EH^{2}(N_{k}(0,\Sigma))<\infty$ for any
invertible $\Sigma\in\mathbb{R}^{k\times k}$. We omit the details.
\end{remark}

\begin{remark}
Ho and Sun \cite{HoSu87} have shown a non-functional version of
Theorem \ref%
{th2} for statistics of the type (\ref{generalmulti}) when the correlation
function $r_{n}$ does not depend on $n$. To the best of our knowledge,
Theorem \ref{th2} is the first central limit theorem for (general)
multipower variation of a row-wise stationary Gaussian process.
\end{remark}

\section{Multipower variation for $\mathcal{BSS}$ processes}\label{sec4}

Armed with the general theorems proved in Section~\ref{sec3}, we are now set to
establish laws of large numbers and central limit results for multipower
variations in the framework of the Brownian semistationary processes. The
regularity conditions invoked are given in a first subsection, while the
next states the theorems, the main parts of the proofs being postponed to
Section~\ref{sec8}; the third subsection discusses the nature of the rather
technical regularity conditions and describes a set of simpler assumptions
that are more amenable to checking.

\subsection{Conditions}\label{sec4.1}

We consider a filtered probability space $(\Omega,\mathcal
{F},(\mathcal{F}%
_{t})_{t\geq0},P)$, assuming the existence thereon of a $\mathcal
{BSS}$\
process, without drift term for the time being, that is,
%
\begin{equation} \label{y}
Y_{t}=\int_{-\infty}^{t}g(t-s)\sigma_{s}W(\mathrm{d}s),
\end{equation}
where $W$ is an $( \mathcal{F}_{t}) $-Brownian measure on
$%
\mathbb{R}$, $\sigma$ is an $( \mathcal{F}_{t})
$-adapted and c%
adlag volatility process and $g\dvtx \mathbb{R}\rightarrow\mathbb{R}$ is a
deterministic continuous memory function with $g(t)=0$ for $t\leq0$
and $%
g\in L^{2}((0,\infty)).$ We also require $\int_{-\infty
}^{t}g^{2}(t-s)\sigma_{s}^{2}\,\mathrm{d}s<\infty$ a.s. to ensure that
$%
Y_{t}<\infty$ a.s. for all $t\geq0$. By an $( \mathcal
{F}_{t}) $%
-Brownian measure we understand a Gaussian stochastic measure such
that, for
any Borelian set $A$ with $E[W(A)^{2}]<\infty, $
\begin{eqnarray*}
W(A)\thicksim N(0,m(A)),
\end{eqnarray*}
where $m$ is the Lebesgue measure, and if $A\subseteq\lbrack t,+\infty),$
then $W(A)$ is independent of $\mathcal{F}_{t}$. Note that $\{
B_{t}:=\int_{a}^{t}W(\mathrm{d}s),t\geq a\} $ is a standard Brownian
motion starting in $a$.

The process $Y$ is assumed to be observed at time points $t_{i}=i/n$, $%
i=1,\ldots,[nt]$. Now, let $G$ be the stationary Gaussian process defined
as
%
\begin{equation}\label{gproc}
G_{t}=\int_{-\infty}^{t}g(t-s)W(\mathrm{d}s).
\end{equation}
This is an important auxiliary object in the study of $\mathcal{BSS}$
processes. Note that $G$ belongs to the type of processes occurring in
(\ref%
{contgaus}), and that the autocorrelation function of $G$ is
%
\begin{equation}\label{autocor}
r(t)=\frac{\int_{0}^{\infty}g(t+u)g(u)\,\mathrm{d}u}{\int
_{0}^{\infty}g^2(u)%
\,\mathrm{d}u}.
\end{equation}
We are interested in the asymptotic behaviour of the functionals
\begin{eqnarray*}
V(Y,p_{1},\ldots,p_{k})_{t}^{n}=\frac{1}{n\tau_{n}^{p_{+}}}%
\sum_{i=1}^{[nt]-k+1}\prod_{j=1}^{k}|\Delta
_{i+j-1}^{n}Y|^{p_{j}} ,\qquad
p_{1},\ldots,p_{k}\geq0 ,
\end{eqnarray*}
where $\Delta_{i}^{n}Y=Y_{{i/n}}-Y_{{(i-1)/n}}$ and $\tau
_{n}^{2}=\bar{R}(1/n)$ with $\bar{R}(t)=E[|G_{s+t}-G_{s}|^{2}]$,
$t\geq0$.
In the following we assume that the function $g$ is continuously
differentiable on $(0,\infty)$, $| g^{\prime}| $ is
non-increasing on $(b,\infty)$ for some $b>0$ and $g^{\prime}\in
L^{2}((\varepsilon,\infty))$ for any $\varepsilon>0$. Moreover, we assume
that for any $t>0,$
%
\begin{equation}\label{finite}
F_{t}=\int_{1}^{\infty}(g^{\prime}(s))^{2}\sigma_{t-s}^{2}\,\mathrm{d}s<\infty
\end{equation}
almost surely.

\begin{remark}
Assumption (\ref{finite}) ensures that the process $Y$ has the same
``smoothness'' as the process $G$ (see Lemma \ref{lemma1} in Section
\ref{sec8}). It is
rather easy to check in practice, because it is implied by the
condition $%
EF_t<\infty$ for $t>0$. Furthermore, if $g$ has bounded support,
assumption (%
\ref{finite}) is trivially fulfilled since $\sigma$ is cadlag.
\end{remark}

\begin{remark}
Let us note again that the process $Y$ is, in general, not a semimartingale.
In particular, this is the case when $g^{\prime}\notin
L^{2}((0,\infty))$%
. For a closer discussion, see \cite{BNSch09}. On the other hand, the
process $Y$ is not of the form (\ref{contgaus}). Thus, we require new
methods to prove the asymptotic results for $V(Y,p_{1},\ldots
,p_{k})_{t}^{n} $. Processes of the form (\ref{y}) are used for modelling
velocity of turbulent flows; see \cite{BNSch07a,BNSch07b,BNSch08a}.
The function $g$, which is used in such models, behaves often as $%
x^{\delta}$ near the origin. Hence, when $\delta\in
(-1/2,1/2)\backslash
\{0\}$, $Y$ is neither a differentiable process nor a semimartingale
(because $g^{\prime}\notin L^{2}((0,\infty))$). This is the primary case
of our interest.
\end{remark}

We define the correlation function of the increments of $G$:
\begin{eqnarray*}
r_{n}(j)=\operatorname{cov}\biggl(\frac{\Delta_{1}^{n}G}{\tau_{n}},\frac
{\Delta
_{1+j}^{n}G}{\tau_{n}}\biggr)=\frac{\bar{R}({(j+1)/n})+\bar
{R}({(j-1)/
n})-2\bar{R}({j/n})}{2\tau_{n}^{2}} ,\qquad j\geq0.
\end{eqnarray*}
Next, we introduce a class of measures that is crucial for our
purposes. We
define (recall that $g(x):=0$ for $x\leq0$)
%
\begin{equation}\label{pimeas}
\pi^{n}(A)=\frac{\int_{A}(g(x-{1/n})-g(x))^{2}\,\mathrm{d}x}{%
\int_{0}^{\infty}(g(x-{1/n})-g(x))^{2}\,\mathrm{d}x} ,\qquad
A\in
\mathcal{B}(\mathbb{R}).
\end{equation}
Note that $\pi^{n}$ is a probability measure on $\mathbb{R}_{+}$.

For the weak law of large numbers we require the following assumptions:

\begin{longlist}[(LLN)]
\item[(LLN)] There exists a sequence $\mathrm{r}(j)$ with
\begin{eqnarray*}
r_{n}^{2}(j)\leq\mathrm{r}(j) ,\qquad\frac{1}{n}\sum
_{j=1}^{n-1}\mathrm{r}%
(j)\rightarrow0.
\end{eqnarray*}
Moreover, it holds that%
%
\begin{equation}\label{measurecond}
\lim_{n\rightarrow\infty}\pi^{n}((\varepsilon,\infty))=0
\end{equation}
for any $\varepsilon>0$.
\end{longlist}

\begin{remark}
\begin{enumerate}[(ii)]
\item[(i)] The first condition of \textup{(LLN)} is adapted from
Theorem %
\ref{th1}. It guarantees the ucp convergence of $V(G,p_{1},\ldots
,p_{k})_{t}^{n}$. The second condition of \textup{(LLN)} says that the whole
mass of the measure $\pi^{n}$ concentrates at $0$. In particular, it is
equivalent to the weak convergence
\begin{eqnarray*}
\pi^{n}\rightarrow\delta_{0} ,
\end{eqnarray*}
where $\delta_{0}$ is the Dirac measure at $0$.

\item[(ii)] Condition (\ref{measurecond}) is absolutely crucial for the
limit theorems given in the next subsection. When this condition is violated
things become more complicated. In particular, it may lead to a different
stochastic limit of $V(Y,p_{1},\ldots,p_{k})_{t}^{n}$ (see the first
example in Section~\ref{sec5}). Intuitively, this can be explained by the observation
that the increments $\Delta_{i}^{n}Y$ contain substantial information about
the volatility (far) outside of the interval $[\frac{i-1}{n},\frac{i}{n}]$
when condition (\ref{measurecond}) does not hold. Thus, in general, we
can not expect the limit described in Theorem \ref{th3} below.
\end{enumerate}
\end{remark}

Now, we introduce the assumptions for the central limit theorem:

\begin{longlist}[(CLT)]
\item[(CLT)] Assumption \textup{(LLN)} holds, and
\begin{eqnarray*}
r_{n}(j)\rightarrow\rho(j) ,\qquad j\geq0,
\end{eqnarray*}
where $\rho(j)$ is the correlation function of $(Q_{i})_{i\geq1}$, as
introduced in (\ref{corrQ}). Furthermore, there exists a sequence
$\mathrm{r}%
(j)$ such that, for any $j,n\geq1$,
\begin{eqnarray*}
r_{n}^{2}(j)\leq\mathrm{r}(j) ,\qquad\sum_{j=1}^{\infty}\mathrm{r}
(j)<\infty,
\end{eqnarray*}
and for some $\gamma\in(0,1]$ we have
%
\begin{equation}\label{smooth}
E[|\sigma_{t}-\sigma_{s}|^{A}]\leq C|t-s|^{A\gamma}
\end{equation}
for any $A>0$. Finally, set $p=\min_{1\leq i\leq k,1\leq j\leq
d}(p_{i}^{j})$%
. Assume that $\gamma(p\wedge1)>\frac{1}{2}$ and \ that there exists a
constant $\lambda<-\frac{1}{p\wedge1}$ such that for any
$\varepsilon
_{n}=\mathrm{O}(n^{-\kappa})$, $\kappa\in(0,1)$, we have
%
\begin{equation}\label{meascond}
\pi^{n}((\varepsilon_{n},\infty))=\mathrm{O}\bigl(n^{\lambda(1-\kappa)}\bigr).
\end{equation}
\end{longlist}

\subsection{Limit theorems}\label{sec4.2}

In this section we present the limit laws of multipower variations of $%
\mathcal{BSS}$ processes, in part widening the scope slightly to allow more
general drift terms. Recall that the (realised) multipower variation of a
process $Y$ of the form (\ref{y}) is defined as
%
\begin{equation}\label{multipower}
V(Y,p_{1},\ldots,p_{k})_{t}^{n}=\frac{1}{n\tau_{n}^{p_{+}}}%
\sum_{i=1}^{[nt]-k+1}\prod_{j=1}^{k}|\Delta
_{i+j-1}^{n}Y|^{p_{j}} ,\qquad
p_{1},\ldots,p_{k}\geq0 ,
\end{equation}
where $\tau_{n}^{2}=\bar{R}( \frac{1}{n}) $ and $%
p_{+}=\sum_{j=1}^{k}p_{j}$. Our first result is the following probability
limit theorem.

\begin{theorem}
\label{th3} Consider a process $Z=Z_{1}+Z_{2}$, where $Z_{2}=Y$ is
given by (%
\ref{y}). Assume that the condition \textup{(LLN)} holds and that
%
\begin{equation}\label{pr1cond}
\frac{1}{n\tau_{n}^{p_{+}}}\sum_{i=1}^{[nt]-k+1}\prod
_{j=1}^{k}|\Delta
_{i+j-1}^{n}Z_{\iota_{j}}|^{p_{j}}\stackrel{P}{\longrightarrow}0,
\end{equation}
where $\iota_{1},\ldots,\iota_{k}\in\{1,2\}$, for any $t>0$ and any
$%
(\iota_{1},\ldots,\iota_{k})\not=(2,\ldots,2)$. Define
\begin{eqnarray*}
\rho_{p_{1},\ldots,p_{k}}^{(n)}=E\biggl[\biggl|{\frac{\Delta
_{1}^{n}G}{\tau
_{n}}}\biggr|^{p_{1}}\cdots\biggl|{\frac{\Delta_{k}^{n}G}{\tau
_{n}}}\biggr|%
^{p_{k}}\biggr].
\end{eqnarray*}
Then we have
%
\begin{equation} \label{stoch2}
V(Z,p_{1},\ldots,p_{k})_{t}^{n}-\rho_{p_{1},\ldots
,p_{k}}^{(n)}\int_{0}^{t}|\sigma_{s}|^{p_{+}}\,\mathrm{d}s\stackrel
{\mathrm{ucp}}{%
\longrightarrow}0.
\end{equation}
\end{theorem}

\begin{pf} See Section \ref{sec8}.
\end{pf}

\begin{remark}
The multipower variation is robust to drift processes $Z_{1}$ that are
smoother than the process $Y$. Assume, for instance, that the process $Z_{1}$
satisfies
\begin{eqnarray*}
E[|Z_{1}(t)-Z_{1}(s)|^{p}]=\mathrm{o}\bigl(\bar{R}^{p/2}(|t-s|)\bigr)
\end{eqnarray*}
for every $p>0$. In this case, condition (\ref{pr1cond}) is obviously
satisfied.
\end{remark}

Next, we demonstrate a joint central limit theorem for a family $%
(V(Z,p_{1}^{j},\ldots,p_{k}^{j})_{t}^{n})_{1\leq j\leq d}$ of multipower
variations.

\begin{theorem}
\label{th4} Consider a process $Z=Z_{1}+Z_{2}$, where $Z_{2}=Y$ is
given by (%
\ref{y}). Assume that the condition \textup{(CLT)} holds
and that
\begin{eqnarray*}
\frac{1}{\sqrt{n}\tau_{n}^{p_{+}^{j}}}\sum_{i=1}^{[nt]-k+1}\prod
_{l=1}^{k}|%
\Delta_{i+l-1}^{n}Z_{i_{l}}|^{p_{l}^{j}}\stackrel{P}{\longrightarrow}0,
\end{eqnarray*}
where $\iota_{1},\ldots,\iota_{k}\in\{1,2\}$, for any $t>0$ and any
$%
(\iota_{1},\ldots,\iota_{k})\not=(2,\ldots,2)$. Then we obtain the
stable convergence
%
\begin{equation}
\sqrt{n}\biggl(V(Z,p_{1}^{j},\ldots,p_{k}^{j})_{t}^{n}-\rho
_{p_{1}^{j},\ldots,p_{k}^{j}}^{(n)}\int_{0}^{t}|\sigma
_{s}|^{p_{+}^{j}}%
\,\mathrm{d}s\biggr)_{1\leq j\leq d}\stackrel{\mathrm{st}}{\longrightarrow}%
\int_{0}^{t}A_{s}^{1/2}\,\mathrm{d}B_{s} ,
\end{equation}
where $B$ is a $d$-dimensional Brownian motion that is defined on an
extension of the filtered probability space $(\Omega,\mathcal
{F},(\mathcal{F%
}_{t})_{t\geq0},P)$ and is independent of $\mathcal{F}$, $A$ is a $%
d\times d$-dimensional process given by
%
\begin{equation}\label{amatrix}
A_{s}^{ij}=\beta_{ij}|\sigma_{s}|^{p_{+}^{i}+p_{+}^{j}} ,\qquad1\leq
i,j\leq d ,
\end{equation}
and the $d\times d$ matrix $\beta$ is defined in (\ref{beta}).
\end{theorem}

\begin{pf} See Section \ref{sec8}.
\end{pf}

\subsection{Discussion of assumptions}\label{sec4.3}

We start our discussion again by considering the auxiliary, centered,
stationary Gaussian, process
\begin{eqnarray*}
G_{t}=\int_{-\infty}^{t}g(t-s)W(\mathrm{d}s).
\end{eqnarray*}
First of all, we want to demonstrate how Theorems 1 and 2 apply for the
multipower variation of the process $G$. In other words, we will give a hint
how to check the conditions of these theorems.

Recall definition (\ref{Rbar}) of the variance function $\bar{R}$ of the
increments of $G$ and note that
\begin{eqnarray*}
\bar{R}(t)=E[|G_{s+t}-G_{s}|^{2}]=\int_{0}^{t}g^{2}(x)\,\mathrm{d}%
x+\int_{0}^{\infty}\bigl(g(t+x)-g(x)\bigr)^{2}\,\mathrm{d}x ,\qquad t\geq0.
\end{eqnarray*}
Clearly, the asymptotic behaviour of the multipower variation of the process
$G$ is fully determined by the behaviour of the function $\bar{R}$
near $0$.
As we deal with a continuous process $G$, it is natural to assume that
$\bar{%
R}(t)$ behaves essentially as $t^{\alpha}$ (for some $\alpha>0$) near $0$
(later on we will formalize this assumption). Since the case where the paths
of $G$ are differentiable (a.s.) is not very interesting for us
(because the
consistency can be deduced by the mean value theorem), we concentrate
on the
region $0<\alpha<2$ (the corresponding $g(t)$ behaving as $t^{(\alpha
-1)/2} $).

Let us introduce a new set of assumptions that correspond to the previous
discussion. These assumptions were proposed by Guyon and Leon in \cite%
{GuyLe89} (those authors considered the case of centered, stationary Gaussian
processes $X$; this relates to the $\mathcal{BSS}$ setting with
$\sigma$
constant) and the same assumptions were used in \cite{BNCP07} and
\cite{BNCPW08}.
\begin{longlist}[(A1)]
\item[(A1)] $\bar{R}(t)=t^{\alpha}L_{0}(t)$ for some
$\alpha\in
(0,2)$ and some positive slowly varying (at 0) function $L_{0}$, which is
continuous on $(0,\infty)$.

\item[(A2)]  $\bar{R}^{\prime\prime}(t)=t^{\alpha-2}L_{2}(t)$
for some slowly varying function $L_{2}$, which is continuous on
$(0,\infty
) $.

\item[(A3)] There exists $b\in(0,1)$ with
\begin{eqnarray*}
K=\limsup_{x\rightarrow0}\sup_{y\in[ x,x^{b}]}\biggl|{\frac
{L_{2}(y)}{%
L_{0}(x)}}\biggr|<\infty.
\end{eqnarray*}
\end{longlist}
Recall that a function $L\dvtx(0,\infty)\rightarrow\mathbb{R}$ is called slowly
varying at 0 when the identity
\begin{eqnarray*}
\lim_{x\searrow0}{\frac{L(tx)}{L(x)}}=1
\end{eqnarray*}
holds for any fixed $t>0$.

Now, note that under assumption \textup{(A1)} we have, for any $j\geq1$,
%
\begin{eqnarray} \label{Converg}
r_{n}(j) &=&\operatorname{cov}\biggl(\frac{\Delta_{1}^{n}G}{\tau_{n}},\frac
{\Delta
_{1+j}^{n}G}{\tau_{n}}\biggr) \nonumber
\\
&=&\frac{\bar{R}({(j+1)/n})+\bar{R}({(j-1)/n})-2\bar
{R}({j/n})%
}{2\bar{R}({1/n})}\rightarrow\rho(j)
\\
&=&\frac
{1}{2}\bigl((j+1)^{\alpha
}-2j^{\alpha}+(j-1)^{\alpha}\bigr) ,\nonumber
\end{eqnarray}
because $L_{0}$ is slowly varying at $0$. It is obvious that $\rho(j)$ is
the correlation function of the discrete-time stationary Gaussian
process $%
Q_{i}=B_{i}^{\alpha/2}-B_{i-1}^{\alpha/2}$, where $B^{\alpha/2}$ is a
fractional Brownian motion with parameter $\alpha/2$.

\begin{remark}
It is easy to see that the convergence of
\[
\frac{\bar{R}(
{(j+1)/n})+\bar{%
R}({(j-1)/n})-2\bar{R}({j/n})}{2\bar{R}({1/n})}
\]
to some
real number for all $j\geq1$ implies that there exits an $s(j)$ such that
\begin{eqnarray*}
\frac{\bar{R}({j/n})}{\bar{R}({1/n})}\rightarrow
s(j).
\end{eqnarray*}
Since the result in Theorem \ref{th4} is independent of the scale of
time we
use, we must have
\begin{eqnarray*}
\frac{\bar{R}(j\Delta)}{\bar{R}(\Delta)}\mathop{%
\rightarrow}\limits_{\Delta\downarrow0}s(j)
\end{eqnarray*}
for any $\Delta$ and then $s(jk)=s(j)s(k)$; consequently
$s(j)=j^{\alpha},$
for a certain $\alpha\in\mathbb{R}$. Moreover, since $(j+1)^{\alpha
}-2j^{\alpha}+(j-1)^{\alpha}$ is a covariance function, we have
$0<\alpha
<2.$ So in the present setting, $(Q_{i})_{i\geq1},$ as defined in
Section~\ref{sec3},
is always a standard fractional Gaussian noise.
\end{remark}

As shown in \cite{GuyLe89,BNCP07} assumptions \textup{(A1)--(A3)}
imply that condition (\ref{pncond}) holds for any $\alpha\in(0,2)$ and
condition (\ref{cltcond2}) holds for any $\alpha\in(0,3/2)$. Hence,
Theorem \ref{th1} holds for all $\alpha\in(0,2)$ while Theorem \ref{th2}
only holds for $\alpha\in(0,3/2)$.

Now, let us see what the conditions \textup{(A1)--(A3)} mean for the memory
function $g$. For simplicity, let us consider functions of the form
%
\begin{equation}\label{functiong}
g(x)=x^{\delta}1_{(0,1]}(x) ,\qquad x>0.
\end{equation}
For such functions we readily obtain assumptions \textup{(A1) and (A2)} with
\begin{eqnarray*}
\alpha=2\delta+1 ,\qquad\delta\in\bigl(-\tfrac{1}{2},0\bigr)\cup\bigl(0,\tfrac{1}{2}\bigr)
\end{eqnarray*}
(the technical assumption \textup{(A3)} has to be checked separately;
for an
example, see Section~\ref{sec5}). Note that for $\delta=0$, for which assumption
\textup{(A2)} does not hold, the process $G$ is a semimartingale and the
multipower variations can be treated as in \cite{BNGJPS06}.

Next, we discuss the assumptions of Section~\ref{sec4.1} for the function $g$
defined in (\ref{functiong}). Recall that condition (\ref{finite}) is
automatically satisfied for functions $g$ with compact support (as in
(\ref%
{functiong})). A straightforward calculation shows that
\begin{eqnarray*}
\pi^{n}((\varepsilon,\infty))=\mathrm{O}((n\varepsilon)^{2\delta-1})
\end{eqnarray*}
for any $\varepsilon>\frac{1}{n}$. Thus, condition (\ref{measurecond})
of \textup{(LLN)} is satisfied (because $2\delta-1<0$) and Theorem
\ref{th3}
is valid for all $\delta\in(-1/2,0)\cup(0,1/2)$.

Finally, we explain how to verify condition (\ref{meascond}) of
\textup{%
{(CLT)}}. Recall that $p=\break \min_{1\leq i\leq k,1\leq j\leq
d}(p_{i}^{j}) $. Let $\varepsilon_{n}=n^{-\kappa}$, $\kappa\in
(0,1)$. We
readily deduce that
\begin{eqnarray*}
\pi^{n}((\varepsilon,\infty))=\mathrm{O}\bigl(n^{\lambda(1-\kappa)}\bigr) ,\qquad
\lambda=%
{2\delta-1}.
\end{eqnarray*}
Thus, condition (\ref{meascond}) is satisfied if
\begin{eqnarray*}
\lambda<-\frac{1}{1\wedge p}.
\end{eqnarray*}
We immediately deduce that Theorem \ref{th4} holds if
\begin{eqnarray*}
&&\hspace*{19pt}p\geq1\mbox{:}\qquad\gamma>\frac{1}{2},\qquad\delta\in\biggl(-\frac{1}{2},0\biggr),
\\
&&\frac{1}{2}<p<1\mbox{:}\qquad\gamma>\frac{1}{2p},\qquad\delta\in\biggl(-\frac
{1}{2},%
\frac{p-1}{2p}\biggr).
\end{eqnarray*}

\begin{remark}
Clearly, we can deal with a larger class of functions $g$ than $%
g(x)=x^{\delta}1_{(0,1]}(x)$. Assume that condition (\ref{finite}) holds.
In the following we consider functions $L_{g}$, $L_{g^{\prime}}$,
which are
continuous on $(0,\infty)$ and slowly varying at 0. We assume the following
conditions:

Assumption: $g\in L^{2}((0,\infty))$ and for some $\delta\in
(-1/2,0)\cup
(0,1/2)$ it holds that:
\begin{enumerate}[(ii)]
\item[(i)] $g(x)=x^{\delta}L_{g}(x)$.

\item[(ii)] $g^{\prime}(x)=x^{\delta-1}L_{g^{\prime}}(x)$ and, for
any $%
\varepsilon>0$, $g^{\prime}\in L^{2}((\varepsilon,\infty))$.
Moreover, $%
| g^{\prime}| $ is non-increasing on $(b,\infty)$ for
some $b>0$%
.
\end{enumerate}

We further assume that the function
\begin{eqnarray*}
\bar{R}(t)=\int_{0}^{t}g^{2}(x)\,\mathrm{d}x+\int_{0}^{\infty}\bigl(g(t+x)-g(x)\bigr)^{2}\,\mathrm{d}x
\end{eqnarray*}
satisfies conditions \textup{(A1)--(A3)} with $\alpha=2\delta+1$.

Under these assumptions we conclude (as for the simple example $%
g(x)=x^{\delta}1_{(0,1]}(x)$) that Theorem \ref{th3} holds for any
$\delta
\in(-1/2,0)\cup(0,1/2)$, and Theorem \ref{th4} holds when further
\begin{eqnarray*}
&&\hspace*{19pt}p\geq1\mbox{:}\qquad\gamma>\frac{1}{2},\qquad\delta\in\biggl(-\frac{1}{2},0\biggr),
\\
&&\frac{1}{2}<p<1\mbox{:}\qquad\gamma>\frac{1}{2p},\qquad\delta\in\biggl(-\frac
{1}{2},%
\frac{p-1}{2p}\biggr).
\end{eqnarray*}
In both cases we have $Q_{i}=B_{i}^{\delta+
{1/2}}-B_{i-1}^{\delta+%
{1/2}}$, $i\geq1$.
\end{remark}

\section{Examples}\label{sec5}

This section discusses two examples of choice of the damping function
$g$ and the associated probabilistic limit behaviour.

As above let $r$\ denote the autocorrelation function of $G=\int
_{-\infty
}^{\cdot}g(\cdot-s)W(\mathrm{d}s)$. Note that assumptions \textup
{(A1)--(A3)%
} could equivalently have been formulated in terms of $1-r$\ rather
than $%
\bar{R}$ \ (since $\bar{R}( t) =2\| g\|
^{2}(1-r(
t) )$).

Suppose first that
\begin{eqnarray*}
g( t) =\mathrm{e}^{-\lambda t}1_{( 0,1) }(
t)
\end{eqnarray*}
with $\lambda>0$. This example (for a detailed discussion, see \cite
{BNSch09}%
) is a non-semimartingale case, and it can be shown that $\pi
^{n}\rightarrow\pi$, with $\pi$ given by
\begin{eqnarray*}
\pi=\frac{1}{1+\mathrm{e}^{-2\lambda}}\delta_{0}+\frac{1}{1+\mathrm{e}^{2\lambda
}}\delta
_{1},
\end{eqnarray*}
where $\delta_{i}$ is the Dirac measure at $i$. Moreover,
\begin{eqnarray*}
V( Y,2) _{t}^{n}\stackrel{P}{\longrightarrow}(
1+\mathrm{e}^{-2\lambda}) ^{-1}\sigma_{t}^{2+}+( 1+\mathrm{e}^{2\lambda
})
^{-1}\sigma_{(-1,t-1]}^{2+},
\end{eqnarray*}
where for any $a<b$
\begin{eqnarray*}
\sigma_{(a,b]}^{2+}=\int_{a}^{b}\sigma_{s}^{2}\,\mathrm{d}s.
\end{eqnarray*}
Thus, in particular, we do not have $V( Y,2)
_{t}^{n}\stackrel{P}{%
\longrightarrow}\sigma_{t}^{2+}$. Note that in this example assumption
\textup{(A2)} is not satisfied.

Our main example is
%
\begin{equation}
g( t) =t^{\nu-1}\mathrm{e}^{-\lambda t}1_{( 0,\infty
) }(
t) \label{g}
\end{equation}
for $\lambda>0$ and with $\nu>\frac{1}{2}$. (So, for $t$ near 0, $g(t)$
behaves as $t^{\delta}$ with $\delta=\nu-1$). \ The autocorrelation
function is given by
%
\begin{equation}\label{modelr}
r( t) =\frac{( 2\lambda) ^{2\nu-1}}{\Gamma
(
2\nu-1) }\mathrm{e}^{-\lambda t}\int_{0}^{\infty}( t+u)
^{\nu
-1}u^{\nu-1}\mathrm{e}^{-2\lambda u}\,\mathrm{d}u.
\end{equation}
It can be proved using properties of Bessel functions (see Sections 5.1--5.3
in \cite{BNCP09b}) that assumptions \textup{(A1)--(A3)} are met
provided that $%
\alpha=2\nu-1\in( 0,2) ,$ that is, $\nu\in(\frac
{1}{2},\frac{3}{2%
}),$ and that $\rho_{p_{1}^{j},\ldots,p_{k}^{j}}^{( n) }$
may be
substituted by $\rho_{p_{1}^{j},\ldots,p_{k}^{j}}$ in the central limit theorem
provided $\nu\in( \frac{1}{2},\frac{5}{4}) $. For
$t\rightarrow
0$ we also have the following asymptotic equivalence
\begin{eqnarray*}
1-r( t) \sim
\cases{2^{-2\nu+1}\dfrac{\Gamma( {3/2}-\nu) }{\Gamma
( \nu
+{1/2}) }( \lambda t) ^{2\nu-1}+\mathrm{O}(
t^{2})
 & \quad\mbox{for}  $\dfrac{1}{2}<\nu<\dfrac{3}{2}$,
 \cr
\dfrac{1}{2}( \lambda t) ^{2}|\log t| & \quad\mbox{for }  $\nu
=\dfrac{3}{2}$,
\cr
\dfrac{1}{4( \nu-{3/2}) }( \lambda t)
^{2}+\mathrm{O}( t^{2\nu-1})  &\quad \mbox{for }  $\dfrac{3}{2}<\nu$.
}
\end{eqnarray*}

\begin{remark}
So for $\frac{3}{2}<\nu\leq2$ the autocorrelation function is twice
differentiable at $0$ and consequently $Y$\ has continuously differentiable
sample paths, while for $\frac{1}{2}<\nu\leq\frac{3}{2}$\ the
sample paths
are Lipschitz of order $\lambda$ for every $0<\lambda<\nu-\frac{1}{2}$
(cf. \cite{CL67}, Section 9.2).
\end{remark}

\section{An application}\label{sec6}

Let us consider the realised variation ratio (RVR) defined for a stochastic
process $X$ as
%
\begin{equation}\label{RVR}
RVR_{t}^{n}:=\frac{{(\uppi/2)}V(X,1,1)_{t}^{n}}{V(X,2,0)_{t}^{n}}.
\end{equation}
The RVR is of interest as a diagnostic tool concerning the nature of
empirical processes.

In particular, it can be used to test the hypothesis that such a
process is
a Brownian semimartingale (with a non-trivial local martingale component)
against the possibility that it is of this type plus a jump process,
see %
\cite{BNSW06,Jac08c} (some related work is discussed in \cite%
{Woe08}). If a jump component is present, then the limit of
$RVR_{t}^{n}$ is
smaller than $1$.

However, in the course of the turbulence project, mentioned earlier, when
calculating the RVR for an extensive high-quality data set from atmospheric
turbulence it turned out that the values of RVR were consistently higher
than $1$. The wish to understand this phenomenon has been a strong
motivation for the theoretical developments described in this paper.
As a consequence of Theorem \ref{th3}, we obtain the following probability
limit result for the realised variation ratio of $\mathcal{BSS}$ processes:
%
\begin{equation}
RVR_{t}^{n}-\psi( r_{n}(1)) \stackrel
{\mathrm{ucp}}{\longrightarrow}0,
\end{equation}
where

\begin{equation}
\psi( \rho) =\sqrt{1-\rho^{2}}+\rho\arcsin\rho,
\end{equation}
which equals $\frac{\uppi}{2} \mathrm{E}\{ | UV|
\} $
of two standard normal variables $U$ and $V$\ with correlation $\rho$.

 Moreover, we have that
\begin{eqnarray*}
\sqrt{n}\bigl(RVR_{t}^{n}-\psi(r_{n}(1))\bigr) &=&\sqrt{n}
\biggl(\frac{{%
(\uppi/2)}V(Y,1,1)_{t}^{n}-\psi(r_{n}(1))\int_{0}^{t}\sigma
_{s}^{2}\,\mathrm{d}%
s}{\int_{0}^{t}\sigma_{s}^{2}\,\mathrm{d}s}\biggr)
\\
&&{}-\sqrt{n}RVR_{t}^{n}\biggl(\frac{V(Y,2,0)_{t}^{n}-\int
_{0}^{t}\sigma_{s}^{2}%
\,\mathrm{d}s}{\int_{0}^{t}\sigma_{s}^{2}\,\mathrm{d}s}\biggr),
\end{eqnarray*}
so, if the parameter $\alpha\in(0,1)$, by applying Theorems \ref
{th3} and %
\ref{th4}, we obtain

\begin{equation}\label{RVRCLT}
\sqrt{n}\bigl(RVR_{t}^{n}-\psi(r_{n}(1))\bigr)\stackrel
{\mathrm{st}}{\longrightarrow}\biggl(%
\frac{\uppi}{2},-\psi(\rho(1))\biggr)\beta^{1/2}\frac{\int_{0}^{t}\sigma
_{s}^{2}%
\,\mathrm{d}B_{s}}{\int_{0}^{t}\sigma_{s}^{2}\,\mathrm{d}s},
\end{equation}
where $\psi$ is as above and the matrix $\beta$ is given in Theorem \ref{th2}.
Specifically, we find $\beta=(\beta_{ij})_{1\leq i,j\leq2}$, where%
\begin{eqnarray*}
\beta_{11} &=&\lim_{n\rightarrow\infty}n \operatorname{var}
(V_{Q}(1,1)_{1}^{n}%
) , \\
\beta_{22} &=&\lim_{n\rightarrow\infty}n \operatorname{var}
(V_{Q}(2,0)_{1}^{n}%
) , \\
\beta_{12} &=&\lim_{n\rightarrow\infty}n \operatorname{cov}(%
V_{Q}(2,0)_{1}^{n},V_{Q}(1,1)_{1}^{n})
\end{eqnarray*}
with $Q$\ as defined in Theorem~\ref{th2}. Thus, we obtain
\begin{eqnarray*}
\beta_{22}=\operatorname{var}(Q_{1}^{2})+2\sum_{k=1}^{\infty}\operatorname{cov}%
(Q_{1}^{2},Q_{1+k}^{2})=2+4\sum_{k=1}^{\infty}\rho^{2}(k).
\end{eqnarray*}
Similarly, we have that
\begin{eqnarray*}
\beta_{12}=2\operatorname{cov}(Q_{1}^{2},|Q_{1}||Q_{2}|)+2\sum_{k=1}^{\infty
}%
\operatorname{cov}(Q_{1}^{2},|Q_{1+k}||Q_{2+k}|).
\end{eqnarray*}
Then, if we write $E[|X_{1}^{2}X_{2}X_{3}|]:=h(\rho_{12},\rho
_{13},\rho
_{23}),$ where $X_{1},X_{2}$ and $X_{3}$ are standard normal with
$\operatorname{cov}%
(X_{i},X_{j})=\rho_{ij}$, we have
\begin{eqnarray*}
\beta_{12}=\biggl(h(1,\rho(1),\rho(1))-\frac{2}{\uppi}\psi(\rho
(1))\biggr)%
+2\sum_{k=1}^{\infty}\biggl(h\bigl(\rho(k),\rho(k+1),\rho(1)\bigr)-\frac
{2}{\uppi}%
\psi(\rho(1))\biggr).
\end{eqnarray*}
To compute the latter, we may use the following formula:
\begin{eqnarray*}
h(\rho_{12},\rho_{13},\rho_{23})=\frac{2}{\uppi}\bigl(\sqrt{1-\rho
_{23}^{2}%
}(1+\rho_{12}^{2}+\rho_{13}^{2})+(\rho_{23}+2\rho_{12}\rho
_{13})\arcsin
(\rho_{23})\bigr);
\end{eqnarray*}
see \cite{Nab52}. For the remaining term we deduce
\begin{eqnarray*}
\beta_{11}=\operatorname{var}(|Q_{1}||Q_{2}|)+2\sum_{k=1}^{\infty}\mbox
{cov}%
(|Q_{1}||Q_{2}|,|Q_{1+k}||Q_{2+k}|).
\end{eqnarray*}
However, while there is no explicit formula available for the latter
expression, it can be computed numerically.

\section{Conclusion and outlook}\label{sec7}

In this paper we have derived convergence in probability and normal
asymptotic limit results for multipower variations of processes $Y$ that,
up to a drift-like term, has the form%
\begin{eqnarray*}
Y_{t}=\int_{-\infty}^{t}g( t-s) \sigma_{s}W(\mathrm{d}s).
\end{eqnarray*}
A key type of example\ has $g( t) $\ behaving as
$t^{\delta}$\
for $t\downarrow0$ and $\delta\in(-\frac{1}{2},\frac
{1}{2})\backslash
\{0\}$. In those instances, $Y$ is not a semimartingale and the limit theory
of multipower variation developed for semimartingales does not suffice to
derive the desired kind of limit results. The basic tool we establish and
apply for this is a normal central limit theorem for triangular arrays of
dependent Gaussian variables. As a case of some special interest for
applications, particularly in turbulence, the central limit behaviour
of the
realised variation ratio, that is, the ratio of bipower variation to quadratic
variation, is briefly discussed. Some specific examples of choice of
$g$ are
also considered.

The turbulence context concerns time-wise observations of
velocities at a single location $x$ in space. More generally, it would
be of
interest to develop the theory of multipower variation corresponding to a
setting where velocities are observed along a curve $\tau$ in space--time.
More specifically, suppose that velocity $Y_{t}( x) $ at position
$x$ and time $t$ is defined by%
\begin{eqnarray*}
Y_{t}( x) =\int_{A+( x,t) }g(
t-s,x-\xi)
\sigma_{s}( \xi) W( \mathrm{d}\xi\,\mathrm
{d}s),
\end{eqnarray*}
where $W$\ denotes white noise, $\sigma_{t}( x) $ is a positive
stationary random field on $\mathbb{R}^{2}$, $g$ is a deterministic damping
function and $A$\ is a subset of space--time involving only points with
negative time coordinates. Then, with the curve $\tau$ parametrized
as $%
\tau( w) =( x( w) ,t( w)
) $,
say, the problem is to study multipower variations of the process $X$\
defined as%
\begin{eqnarray*}
X_{w}=\int_{A+\tau( w) }g\bigl( t( w)
-s,x(
w) -\xi\bigr) \sigma_{s}( \xi) W(
\mathrm{d}\xi
\,\mathrm{d}s) .
\end{eqnarray*}
Among the questions that this raises is that of proper definition of
filtrations.

In another direction it would be of interest to extend the results of
this paper to power and multipower variations of higher-order differences
of $Y$. In particular, this might yield normal central limit theorems for
the whole range of values of $\delta$ and it could also lead to more
robustness against drift processes. For some recent work on quadratic
variation of higher-order differences, see \cite{Beg07a,Beg07b}
and references given there.

\section{Proofs}\label{sec8}

As our proofs are rather long and technical, let us briefly outline the
scheme. The proofs basically consist of four main steps:
\begin{longlist}[(iii)]
\item[(i)] First, we provide a Wiener chaos decomposition for the functional $%
V(p_{1}^{j},\ldots,p_{k}^{j})_{t}^{n}$, $j=1,\ldots,d$, that appears in
Theorems \ref{th1} and \ref{th2} (see Section 8.3).

\item[(ii)] In a second step, we prove the weak law of large numbers and the central
limit theorem for the normalized version of $(V(p_{1}^{j},\ldots
,p_{k}^{j})_{t}^{n})_{1\leq j\leq d}$, using its Wiener chaos decomposition
and recent techniques from Malliavin calculus derived in \cite
{NuPec05,NuO-L08,PecTu05} (Sections~\ref{sec8.4} and \ref{sec8.5}).

\item[(iii)] In order to prove Theorems \ref{th3} and \ref{th4}, we first
show that
the true increments $\Delta_{i}^{n}Y$ can be replaced by the quantity $
\sigma_{{(i-1)/n}}\Delta_{i}^{n}G$ without changing the asymptotic
limits (see (\ref{conv1}) and (\ref{convclt})). For this step, the conditions
(\ref{measurecond}) and (\ref{meascond}) on the measure $\pi^{n}$ are
absolutely crucial.

\item[(iv)] In the last step, we apply the following blocking technique: We divide
the interval $[0,t]$ into big sub-blocks (whose lengths still converge
to $0$%
) and freeze the volatility process $\sigma$ at the beginning of each
sub-block. Then Theorem~\ref{th4} (respectively, Theorem \ref{th3})
follows from
Theorem~\ref{th2} (respectively, Theorem \ref{th1}) applied to the
family of
functionals $(V(G,p_{1}^{j},\ldots,p_{k}^{j})_{t}^{n})_{1\leq j\leq
d}$ and
the properties of stable convergence.
\end{longlist}
Below, all positive constants (which do not depend on $n$) are denoted
by $C$%
, although they might change from line to line.

\subsection{Some elements of Malliavin calculus}\label{sec8.1}

Before we proceed with the proofs of the main results, we review the basic
concepts of the Wiener chaos expansion. Consider a complete probability
space $(\Omega,\mathcal{F},P)$ and a subspace $\mathcal
{H}_{1}\mathcal{\ }$%
of $L^{2}(\Omega,\mathcal{F},P)$ whose elements are zero-mean Gaussian
random variables$.$ Let $\mathbb{H}$ be a separable Hilbert space with scalar
product denoted by $\langle\cdot,\cdot\rangle_{\mathbb{H}}$ and
norm $%
\|\cdot\|_{\mathbb{H}}.$ We will assume that there is an isometry
\begin{eqnarray*}
 W\dvtx \mathbb{H} &\rightarrow&\mathcal{H}_{1}, \\
h &\mapsto&W(h)
\end{eqnarray*}
in the sense that
\begin{eqnarray*}
E[W(h_{1})W(h_{2})]=\langle h_{1},h_{2}\rangle_{\mathbb{H}}.
\end{eqnarray*}
It is easy to see that this map has to be linear.

For any $m\geq2$, we denote by $\mathcal{H}_{m\mbox{ }}$the $m$th Wiener
chaos, that is, the closed subspace of $L^{2}(\Omega,\mathcal{F},P)$
generated by the random variables $H_{m}(X)$, where $X\in\mathcal
{H}_{1}$, $%
E[X^{2}]=1$ and $H_{m}$ is the $m$th Hermite polynomial, that is, $H_{0}(x)=1$
and $H_{m}(x)=(-1)^{m}\mathrm{e}^{{x^{2}/2}}{\frac
{\mathrm{d}^{m}}{\mathrm{d}x^{m}}}(\mathrm{e}^{-{{%
x^{2}/2}}}).$

We denote by
\begin{eqnarray*}
 I_{m}\dvtx \mathbb{H}^{\odot m}\rightarrow\mathcal
{H}_{m\mbox{ }}
\end{eqnarray*}
the isometry between the symmetric tensor product $\mathbb{H}^{\odot m}$,
equipped with the norm $\sqrt{m!}\| \cdot\| _{\mathbb
{H}^{\otimes
m}} $, and the $m$th chaos $\mathcal{H}_{m}$; see Section
1.1.1 in %
\cite{Nu06} for its definition. For $h\in\mathbb{H}^{\otimes m},$ we
set $%
I_{m}(h):=I_{m}(\tilde{h})$, where $\tilde{h}$ is the symmetrization
of $h$.
For any $g\in\mathbb{H}^{\otimes m}$, $h\in\mathbb{H}^{\otimes n}$,
$n,m\geq0$,
it holds that
\begin{eqnarray*}
E[I_{m}(g)I_{n}(h)]=\delta_{mn}m!\langle\tilde{g},\tilde{h}\rangle
_{\mathbb{H}^{\otimes m}}.
\end{eqnarray*}
For any $h=h_{1}\otimes\cdots\otimes h_{m}$ and $g=g_{1}\otimes
\cdots
\otimes g_{m}\in\mathbb{H}^{\otimes m}$, we define the $p$th contraction
of $h$ and $g$, denoted by $h\otimes_{p}g$, as the element of $%
\mathbb{H}^{\otimes2(m-p)}$ given by
\begin{eqnarray*}
h\otimes_{p}g=\langle h_{1},g_{1}\rangle_{\mathbb{H}}\cdots\langle
h_{p},g_{p}\rangle_{\mathbb{H}}h_{p+1}\otimes\cdots\otimes
h_{m}\otimes
g_{p+1}\otimes\cdots\otimes g_{m}.
\end{eqnarray*}
This definition can be extended by linearity to any element of $%
\mathbb{H}^{\otimes m}$.

Now, let $\mathcal{G}$ be the $\sigma$-field generated by the random
variables $\{W(h)| h\in\mathbb{H}\}$. Any square-integrable random
variable $%
F\in L^{2}(\Omega,\mathcal{G},P)$ has a unique chaos decomposition
\begin{eqnarray*}
F=\sum_{m=0}^{\infty}I_{m}(h_{m}) ,
\end{eqnarray*}
where $h_{m}\in \mathbb{H}^{\odot m}$ (see \cite{Nu06} for more details).

Finally, we adapt the theory of Wiener chaos expansion to the set up of
Section~\ref{sec3}. Let $\mathcal{G}$ be the $\sigma$-field generated by the random
variables $(X_{i,n})_{n\geq1,1\leq i\leq\lbrack nt]}$ and $\mathcal{H}_{1}$
be the first Wiener chaos associated with $(X_{i,n})_{n\geq1,1\leq
i\leq
\lbrack nt]}$, that is, the closed subspace of $L^{2}(\Omega,\mathcal{G},P)$
generated by the random variables $(X_{i,n})_{n\geq1,1\leq i\leq
\lbrack
nt]}$. Notice that $\mathcal{H}_{1}$ can be seen as a separable Hilbert
space with a scalar product induced by the covariance function of the
process $(X_{i,n})_{n\geq1,1\leq i\leq\lbrack nt]}$. This means we can
apply the above theory of Wiener chaos expansion with the canonical Hilbert
space $\mathbb{H}=\mathcal{H}_{1}$. Denote by $\mathcal{H}_{m}$ the
$m$th Wiener
chaos associated with the triangular array $(X_{i,n})_{n\geq1,1\leq
i\leq
\lbrack nt]}$ and by $I_{m}$ the corresponding linear isometry between the
symmetric tensor product $\mathcal{H}_{1}^{\odot m}$ (equipped with
the norm
$\sqrt{m!}\| \cdot\| _{\mathcal{H}_{1}^{\otimes m}}$)
and the $m$%
th Wiener chaos.

\subsection{Preliminary results}\label{sec8.2}

First of all, let us note that w.l.o.g. the volatility process $\sigma
$ can
be assumed to be bounded on compact intervals because $\sigma$ is cadl%
ag and it is integrated with respect to $W$. This follows by a standard
localization procedure presented in \cite{BNGJPS06}. Furthermore, the
process $F_{t}$, defined by (\ref{finite}), is continuous, because
$\sigma$
is cadlag. Hence, $F_{t}$ is locally bounded and can be assumed to
be bounded on compact intervals w.l.o.g. by the same localization procedure.

Next we establish three lemmas.

\begin{lemma}
\label{lemma1} Under assumption (\ref{finite}), it holds that
%
\begin{equation} \label{incrass}
E[|\Delta_i^n Y|^p]\leq C_p \tau_n^p , \qquad i=0,\ldots, [nt],
\end{equation}
for all $p>0$.
\end{lemma}

\begin{pf} Recall that $g^{\prime}$ is
non-increasing on $(b,\infty)$ for some $b>0$. Assume w.l.o.g. that $b>1$.
Observe the decomposition
\begin{eqnarray*}
\Delta_{i}^{n}Y=\int_{{(i-1)/n}}^{{i/n}}g\biggl(\frac
{i}{n}-s\biggr)\sigma
_{s}W(\mathrm{d}s)+\int_{-\infty}^{{(i-1)/n}}\biggl(g\biggl(\frac
{i}{n}-s\biggr)-g\biggl(%
\frac{i-1}{n}-s\biggr)\biggr)\sigma_{s}W(\mathrm{d}s).
\end{eqnarray*}
Since $\sigma$ is bounded on compact intervals, we deduce by Burkholder's
inequality
\begin{eqnarray*}
E[|\Delta_{i}^{n}Y|^{p}]\leq C_{p}\biggl( \tau_{n}^{p}+E\biggl(%
\int_{0}^{\infty}\biggl(g\biggl(\frac{1}{n}+s\biggr)-g(s)\biggr)^{2}\sigma_{
{(i-1)/n}%
-s}^{2}\,\mathrm{d}s\biggr)^{p/2}\biggr).
\end{eqnarray*}
We immediately obtain the estimates
\begin{eqnarray*}
\int_{0}^{1}\biggl(g\biggl(\frac{1}{n}+s\biggr)-g(s)\biggr)^{2}\sigma_{
{(i-1)/n}-s}^{2}%
\,\mathrm{d}s&\leq& C\tau_{n}^{2} ,
\\
\int_{1}^{b}\biggl(g\biggl(\frac{1}{n}+s\biggr)-g(s)\biggr)^{2}\sigma_{
{(i-1)/n}-s}^{2}%
\,\mathrm{d}s&\leq&\frac{C}{n^{2}} ,
\end{eqnarray*}
because $g^{\prime}$ is continuous on $(0,\infty)$ and $\sigma$ is
bounded on compact intervals. On the other hand, since $g^{\prime}$ is
non-increasing on $(b,\infty)$, we get
\begin{eqnarray*}
\int_{b}^{\infty}\biggl(g\biggl(\frac{1}{n}+s\biggr)-g(s)\biggr)^{2}\sigma_{
{(i-1)/n}%
-s}^{2}\,\mathrm{d}s\leq\frac{F_{{(i-1)/n}}}{n^{2}}.
\end{eqnarray*}
The boundedness of the process $F$ implies (\ref{incrass}).
\end{pf}

Next, for any stochastic process $f$ and any $s>0$, we define the (possibly
infinite) measure (recall that $g(x):=0$ for $x\leq0$)
%
\begin{equation} \label{pimeasure}
\pi_{f,s}^{n}(A)=\frac{E\int_{A}(g(x-
{1/n})-g(x))^{2}f_{s-x}^{2}%
\,\mathrm{d}x}{\int_{0}^{\infty}(g(x-{1/n})-g(x))^{2}\,\mathrm
{d}x}%
 ,\qquad A\in\mathcal{B}(\mathbb{R}).
\end{equation}

\begin{lemma}
\label{lemma2} Under assumption (\ref{finite}), it holds that
%
\begin{equation}\label{piappr}
\sup_{s\in\lbrack0,t]}\pi_{\sigma,s}^{n}((\varepsilon,\infty
))\leq
C\pi^{n}((\varepsilon,\infty))
\end{equation}
for any $\varepsilon>0$, where $\pi^{n}$ is given by (\ref{pimeas}).
\end{lemma}

\begin{pf} Recall again that $|
g^{\prime
}| $ is non-increasing on $(b,\infty)$ for some $b>0$, and assume
w.l.o.g. that $b>\varepsilon$. Since the processes $\sigma$ and $F$ are
bounded we deduce exactly as in the previous proof that
\begin{eqnarray*}
&&\int_{\varepsilon}^{\infty}\biggl(g\biggl(x-\frac{1}{n}\biggr)-g(x)\biggr)^{2}\sigma
_{s-x}^{2}%
\,\mathrm{d}x
\\
&&\quad =\int_{\varepsilon}^{b}\biggl(g\biggl(x-\frac
{1}{n}\biggr)-g(x)\biggr)^{2}\sigma
_{s-x}^{2}\,\mathrm{d}x+\int_{b}^{\infty}\biggl(g\biggl(x-\frac
{1}{n}\biggr)-g(x)\biggr)^{2}\sigma
_{s-x}^{2}\,\mathrm{d}x
 \\
&&\quad\leq C\biggl(\int_{\varepsilon}^{\infty}\biggl(g\biggl(x-\frac
{1}{n}\biggr)-g(x)\biggr)^{2}\,\mathrm{%
d}x+n^{-2}\biggr).
\end{eqnarray*}
This completes the proof of Lemma \ref{lemma2}.
\end{pf}

Finally, we present the following technical lemma.

\begin{lemma}
\label{lemma3} Under the assumption \textup{{(CLT)}}, there
exists a
number $l\geq1$ and positive sequences $\varepsilon
_{n}^{(j)}\rightarrow0$%
, $j=1,\ldots,l$, such that $0<\varepsilon_{n}^{(1)}<\cdots
<\varepsilon
_{n}^{(l)}$ and
%
\begin{eqnarray}\label{meascond1}
&&\varepsilon_{n}^{(1)}=\mathrm{o}\bigl(n^{-{1/(2\gamma(p\wedge1))}}\bigr) ,\qquad
\pi
^{n}\bigl(\bigl(\varepsilon_{n}^{(l)},\infty\bigr)\bigr)=\mathrm{o}\bigl(n^{-{1/(p\wedge1)}}\bigr),
\\\label{meascond2}
&&\bigl(\varepsilon_{n}^{(j+1)}\bigr)^{2\gamma}\pi^{n}\bigl(\bigl(\varepsilon
_{n}^{(j)},\infty\bigr)\bigr)=\mathrm{o}\bigl(n^{-{1/(p\wedge1)}}\bigr) ,\qquad j=1,\ldots,l-1 ,
\end{eqnarray}
where $p=\min_{1\leq i\leq k,1\leq j\leq d}(p_{i}^{j})$.
\end{lemma}

\begin{pf} Assume first that $p\geq1$. Recall
that $\gamma>1/2$. Set $\varepsilon_{n}^{(j)}=n^{-\kappa_{j}}$, $%
j=1,\ldots,l$, with $1>\kappa_{1}>\cdots>\kappa_{l}>0$. The
condition $%
\pi^{n}((\varepsilon_{n}^{(j)},\infty))=\mathrm{O}(n^{\lambda(1-\kappa_{j})})$
for some $\lambda<-1$, presented in (\ref{meascond}), implies that
conditions (\ref{meascond1}) and (\ref{meascond2}) are satisfied if
we find $%
1>\kappa_{1}>\cdots>\kappa_{l}>0$ such that
\begin{eqnarray*}
\kappa_{1}&>&\frac{1}{2\gamma} , \qquad
\kappa_{l}<1+\frac{1}{\lambda} ,
\\
(1+\lambda)-\kappa_{j}\lambda-2\kappa_{j+1}\gamma&<&0 ,\qquad
1\leq
j\leq l-1.
\end{eqnarray*}
From the first and the last inequality, we deduce by induction that
\begin{eqnarray*}
\frac{1}{2\gamma}<\kappa_{1}<\frac{1+\lambda}{\lambda}\sum
_{i=0}^{l-1}%
\biggl(-\frac{2\gamma}{\lambda}\biggr)^{i}+\biggl(-\frac{2\gamma
}{\lambda}\biggr)%
^{l}\kappa_{l}
\end{eqnarray*}
must hold.

When $2\gamma\geq-\lambda, $ the term on the right-hand side
converges to $%
\infty$ as $l\rightarrow\infty$. In that case it is easy to find
constants $1>\kappa_{1}>\cdots>\kappa_{l}>0$ such that (\ref{meascond1})
and (\ref{meascond2}) are satisfied.

When $2\gamma<-\lambda,$ the limit of $\frac{1+\lambda}{\lambda}%
\sum_{i=0}^{l-1}(-\frac{2\gamma}{\lambda})^{i}$ is $\frac
{%
1+\lambda}{\lambda+2\gamma}$ (as $l\rightarrow\infty$) and the
restriction on $\kappa_{1}$ becomes
\begin{eqnarray*}
\frac{1}{2\gamma}<\kappa_{1}<\frac{1+\lambda}{\lambda+2\gamma}.
\end{eqnarray*}
Notice that $\frac{1}{2\gamma}<\frac{1+\lambda}{\lambda+2\gamma}$
because $\gamma>1/2$. The existence of the positive powers $\kappa
_{j}$, $%
j=2,\ldots,l$ that satisfy the original inequality follows by an induction
argument.

Assume now that $p<1$. Recall that $\gamma$ must satisfy
\begin{eqnarray*}
\gamma>\frac{1}{2p}
\end{eqnarray*}
and $\lambda<-\frac{1}{p}$. Again the conditions (\ref{meascond1})
and (\ref%
{meascond2}) are satisfied if we find $1>\kappa_{1}>\cdots>\kappa_{l}>0$
such that
\begin{eqnarray*}
\kappa_{1}&>&\frac{1}{2\gamma p} , \qquad
\kappa_{l}<\frac{1+\lambda p}{\lambda p} ,
\\
\biggl(\frac{1}{p}+\lambda\biggr) -\kappa_{j}\lambda-2\kappa_{j+1}\gamma
&<&0,\qquad
1\leq j\leq l-1.
\end{eqnarray*}
Notice that the second inequality has solutions because $\lambda
<-\frac{1}{p}
$. Moreover, we deduce as above that the inequality
\begin{eqnarray*}
\frac{1}{2\gamma p}<\kappa_{1}<\frac{{1/p}+\lambda}{\lambda}
\sum_{i=0}^{l-1}\biggl(-\frac{2\gamma}{\lambda}\biggr)^{i}+
\biggl(-\frac{2\gamma
}{\lambda}\biggr)^{l}\kappa_{l}
\end{eqnarray*}
must hold. Again the more complicated case is $2\gamma<-\lambda$. By
letting $l\rightarrow\infty, $ the restriction on $\kappa_{1}$ becomes
\begin{eqnarray*}
\frac{1}{2\gamma p}<\kappa_{1}<\frac{{1/p}+\lambda}{\lambda
+2\gamma
}.
\end{eqnarray*}
Note that $\frac{1}{2\gamma p}<\frac{{1/p}+\lambda}{\lambda
+2\gamma}
$ because $\gamma>\frac{1}{2p}$. As before, the existence of the positive
powers $\kappa_{j}$, $j=2,\ldots,l$ that satisfy the original inequality
follows by an induction argument.
\end{pf}

\subsection{Some notation}

In this subsection we introduce various notation connected to the Wiener
chaos expansion for the functionals $V(p_{1}^{j},\ldots
,p_{k}^{j})_{t}^{n}$%
, $j=1,\ldots,d$, and present some first convergence results.

Recall that the covariance matrix of $(X_{i,n}, \ldots, X_{i+l,n})$ is
invertible for any $l\geq1$ and $n\geq1$. Let $X_i^n(1), \ldots, X_i^n(k)$
be an i.i.d. $N(0,1)$ sequence that spans the same linear space as $X_{i,n},
\ldots, X_{i+k-1,n}$ (such a sequence can be constructed by the Gram--Schmidt
method). Thus, it has the representation
%
\begin{equation} \label{orthrep}
X_i^n(j) = \sum_{l=1}^k a_{lj}^n X_{i+l-1,n} , \qquad j=1, \ldots, k ,
\end{equation}
for some real numbers $a_{lj}^n$. Note that
\begin{eqnarray*}
| a_{lj}^n |\leq C
\end{eqnarray*}
for all $l,j,n$, because $E[X_{i,n}^2]=1$ for all $i,n$.

For any $1\leq j\leq d$, we obtain the Wiener chaos representation
%
\begin{equation}\label{chaosrep}
V(p_{1}^{j},\ldots,p_{k}^{j})_{t}^{n}-\rho_{p_{1}^{j},\ldots
,p_{k}^{j}}^{(n)}t=\sum_{m=2}^{\infty}I_{m}\Biggl(\frac{1}{n}%
\sum_{i=1}^{[nt]}f_{m,j}^{n}(i)\Biggr)+\mathrm{O}_{p}(n^{-1}) ,
\end{equation}
where the $f_{m,j}^{n}(i)\in\mathbb{H}^{\odot m}$ are given by
%
\begin{equation}\label{fform}
f_{m,j}^{n}(i)=\sum_{k_{l}\in\{1,\ldots,k\}}c_{k_{1},\ldots
,k_{m}}^{n}(j)X_{i}^{n}(k_{1})\otimes\cdots\otimes X_{i}^{n}(k_{m})
\end{equation}
for some coefficients $c_{k_{1},\ldots,k_{m}}^{n}(j)$. We set
%
\begin{equation}\label{cform}
c_{m}^{n}(j)=\|f_{m,j}^{n}(i)\|_{\mathbb{H}^{\otimes m}}^{2}=\sum
_{k_{l}\in
\{1,\ldots,k\}}|c_{k_{1},\ldots,k_{m}}^{n}(j)|^{2}.
\end{equation}
Note that
%
\begin{equation}\label{varform}
\operatorname{var}(|X_{i,n}|^{p_{1}^{j}}\cdots
|X_{i+k-1,n}|^{p_{k}^{j}})%
=\sum_{m=2}^{\infty}m!c_{m}^{n}(j)<C
\end{equation}
for all $n,j$, because $E[X_{i,n}^{2}]=1$ for all $i,n$. Finally, when $
f_{m,j}^{n}(i)$, $c_{k_{1},\ldots,k_{m}}^{n}(j)$ and $c_{m}^{n}(j)$
correspond to some particular choice of powers $p_{1},\ldots,p_{k}$,
we use
the notation $f_{m}^{n}(i)$, $c_{k_{1},\ldots,k_{m}}^{n}$ and $c_{m}^{n}$.

Now assume that the assumptions (\ref{cltcond1}) and (\ref{cltcond2}) of
Theorem \ref{th2} hold. Since $a_{lj}^{n}$ in (\ref{orthrep}) is a
continuous function of $r(1),\ldots,r_{n}(k-1)$ and the Gaussian
process $Q$
is non-degenerate, we have that
\begin{eqnarray*}
a_{lj}^{n}\rightarrow a_{lj} ,
\end{eqnarray*}
and the sequence $Q_{i}(1),\ldots,Q_{i}(k)$ given by
%
\begin{equation}\label{qorthrep}
Q_{i}(j)=\sum_{l=1}^{k}a_{lj}Q_{i+l-1} ,\qquad j=1,\ldots,k ,
\end{equation}
is an i.i.d. $N(0,1)$ sequence. Now, let us associate $f_{m,j}(i)$, $%
c_{k_{1},\ldots,k_{m}}(j)$ and $c_{m}(j)$ with the functional $%
V_{Q}(p_{1}^{j},\ldots,p_{k}^{j})_{t}^{n}-\rho_{p_{1}^{j},\ldots
,p_{k}^{j}}t$, where
\begin{eqnarray*}
\rho_{p_{1}^{j},\ldots,p_{k}^{j}}=E[|Q_{1}|^{p_{1}^{j}}\cdots
|Q_{k}|^{p_{k}^{j}}] ,
\end{eqnarray*}
by (\ref{chaosrep}), (\ref{fform}) and (\ref{cform}). By a repeated
application of the multiplication formula (see \cite{Nu06}), we know
that $%
c_{k_{1},\ldots,k_{m}}^{n}(j)$ is a continuous function of
$r_{n}(1),\ldots
,r_{n}(k-1)$. Since $r_{n}(j)\rightarrow\rho(j)$ we obtain
%
\begin{eqnarray}\label{converg1}
c_{k_{1},\ldots,k_{m}}^{n}(j)&\rightarrow& c_{k_{1},\ldots
,k_{m}}(j) ,\qquad
c_{m}^{n}(j)\rightarrow c_{m}(j) ,
\\ \label{converg2}
\langle f_{m,j_{1}}^{n}(i),f_{m,j_{2}}^{n}(i+l)\rangle_{\mathbb
{H}^{\otimes
m}}&\rightarrow&\langle f_{m,j_{1}}(i),f_{m,j_{2}}(i+l)\rangle
_{\mathbb{H}^{\otimes m}},
\end{eqnarray}
\begin{eqnarray}\label{converg3}
&&\operatorname{cov}(|X_{i,n}|^{p_{1}^{j_{1}}}\cdots
|X_{i+k-1,n}|^{p_{k}^{j_{1}}},|X_{i,n}|^{p_{1}^{j_{2}}}\cdots
|X_{i+k-1,n}|^{p_{k}^{j_{2}}})\nonumber
\\
&&\quad=\sum_{m=2}^{\infty}m!\langle
f_{m,j_{1}}^{n}(1),f_{m,j_{2}}^{n}(1)\rangle_{\mathbb{H}^{\otimes m}}
\nonumber
\\[-8pt]\\[-8pt]
&&\qquad\rightarrow\quad\hspace*{-2pt}\operatorname{cov}(|Q_{i}|^{p_{1}^{j_{1}}}\cdots
|Q_{i+k-1}|^{p_{k}^{j_{1}}},|X_{i,n}|^{p_{1}^{j_{2}}}\cdots
|Q_{i+k-1}|^{p_{k}^{j_{2}}}) \nonumber
\\
&&\qquad\qquad\hspace*{1pt}\quad =\sum_{m=2}^{\infty}m!\langle f_{m,j_{1}}(1),f_{m,j_{2}}(1)\rangle
_{\mathbb{H}^{\otimes m}}.\nonumber
\end{eqnarray}

\subsection{Proof of Theorems 1 and 3}\label{sec8.4}

\begin{pf*}{Proof of Theorems \ref{th1}} Since $V(p_1, \ldots, p_k)_t^n$ is
increasing in $t$ and the process $\rho_{p_1, \ldots, p_{k}}^{(n)} t$ is
continuous in $t$, it is sufficient to prove $V(p_1, \ldots, p_k)_t^n -
\rho_{p_1, \ldots, p_{k}}^{(n)} t \stackrel{P}{\longrightarrow} 0$
for a
fixed $t>0$.

Note that
\begin{eqnarray}\label{bracketest}
|\langle f_{m}^{n}(1),f_{m}^{n}(1+l)\rangle_{\mathbb{H}^{\otimes
m}}|&\leq&
c_{m}^{n} ,\nonumber
\\[-8pt]\\[-8pt]
|\langle f_{m}^{n}(1),f_{m}^{n}(1+l)\rangle
_{\mathbb{H}^{\otimes m}}|&\leq& c_{m}^{n}C^{m}\bigl(|r_{n}(l)|^{m}+\cdots
+|r_{n}(l-k+1)|^{m}\bigr),\nonumber
\end{eqnarray}
where the bounds are not comparable. Now, due to assumption (\ref
{pncond}), $%
\mathrm{r}(j)\rightarrow0$ as $j\rightarrow\infty$. Thus, there
exists an $%
H$ such that $|C\mathrm{r}^{1/2}(j-k+1)|<1$ for $j\geq H$ (for any
fixed $C$%
). By (\ref{bracketest}) we have (for any $m\geq2$)
\begin{eqnarray}\label{impest}
\sum_{l=1}^{n-1}|\langle f_{m}^{n}(1),f_{m}^{n}(1+l)\rangle
_{\mathbb{H}^{\otimes m}}|&\leq& C\Biggl(Hc_{m}^{n}+\sum
_{l=H}^{n-1}|\langle
f_{m}^{n}(1),f_{m}^{n}(1+l)\rangle_{\mathbb{H}^{\otimes m}}|\Biggr)\nonumber
\\[-8pt]\\[-8pt]
&\leq& Cc_{m}^{n}\Biggl(H+\sum_{l=H}^{n-1}(C|r_{n}(l)|)^{2}\Biggr)\leq
Cc_{m}^{n}\sum_{l=1}^{n-1}\mathrm{r}(l).\nonumber
\end{eqnarray}
Hence,
\begin{eqnarray*}
\operatorname{var}(V(p_{1},\ldots,p_{k})_{t}^{n})\leq\frac{C}{n}\sum
_{m=2}^{\infty
}m!c_{m}^{n}\Biggl(1+\sum_{l=1}^{n-1}\mathrm{r}(l)\Biggr).
\end{eqnarray*}
The latter converges to 0 due to (\ref{varform}) and assumption (\ref{pncond}).
\end{pf*}

\begin{pf*}{Proof of Theorem \ref{th3}} Assume first that $Z_{1}=0$. Recall
that
%
\begin{equation}\label{incrin}
E[|\Delta_{i}^{n}Y|^{q}]\leq C\tau_{n}^{q} ,\qquad E[|\Delta
_{i}^{n}G|^{q}]\leq C\tau_{n}^{q}
\end{equation}
for any $q\geq0$, due to Lemma \ref{lemma1}.

We assume for simplicity that $k=1$, $p_{1}=p$. The general case can be
proved in a similar manner by (\ref{incrin}) and an application of the
H\"{o}%
lder inequality.

Since $V(Y,p)_{t}^{n}$ is increasing in $t$ and the limit process is
continuous in $t$, it suffices to prove the pointwise convergence $%
V(Y,p)_{t}^{n}\stackrel{P}{\longrightarrow}\mu_{p}\int
_{0}^{t}|\sigma
_{s}|^{p}\,\mathrm{d}s$. For any $l\leq n$, we have
\begin{eqnarray*}
V(Y,p)_{t}^{n}-\mu_{p}\int_{0}^{t}|\sigma_{s}|^{p}\,\mathrm{d}s=\frac
{1}{%
n\tau_{n}^{p}}\sum_{i=1}^{[nt]}\bigl(|\Delta_{i}^{n}Y|^{p}-\big|\sigma
_{{(i-1)/n}}\Delta_{i}^{n}G\big|^{p}\bigr)+R_{t}^{n,l},
\end{eqnarray*}
where
\begin{eqnarray*}
R_{t}^{n,l} &=&{\frac{1}{n\tau_{n}^{p}}}\Biggl(\sum
_{i=1}^{[nt]}|\sigma_{%
{(i-1)/n}}\Delta_{i}^{n}G\big|^{p}-\sum_{j=1}^{[lt]}\bigl|\sigma_{
{(j-1)/l}%
}\bigr|^{p}\sum_{i\in I_{l}(j)}|\Delta_{i}^{n}G|^{p}\Biggr)
\\
&&{}+{\frac{1}{n\tau_{n}^{p}}}\sum_{j=1}^{[lt]}\bigl|\sigma_{{(j-1)/l}
}\bigr|^{p}\sum_{i\in I_{l}(j)}|\Delta_{i}^{n}G|^{p}-\mu
_{p}l^{-1}\sum_{j=1}^{[lt]}\bigl|\sigma_{{(j-1)/l}}\bigr|^{p}
\\
&&{}+\mu_{p}\Biggl(l^{-1}\sum_{j=1}^{[lt]}\bigl|\sigma_{{(j-1)/l}%
}\bigr|^{p}-\int_{0}^{t}|\sigma_{s}|^{p}\, \mathrm{d}s\Biggr)
\end{eqnarray*}
and
\begin{eqnarray*}
I_{l}(j)=\biggl\{i\Big| {\frac{i}{n}}\in\biggl({\frac{j-1}{l}},{\frac
{j}{l}}\biggr]%
\biggr\} ,\qquad j\geq1.
\end{eqnarray*}
The assumption \textup{(LLN)} implies that $V(G,p)_{t}^{n}\stackrel
{\mathrm{ucp}}{%
\longrightarrow}\mu_{p}t$. Since $\sigma$ is cadlag and bounded
on compact intervals, we deduce that
\begin{eqnarray*}
\lim_{l\rightarrow\infty}\lim_{n\rightarrow\infty
}P(|R_{t}^{n,l}|>\epsilon)=0
\end{eqnarray*}
for any $\epsilon>0$. Hence, we are left to prove that
\begin{eqnarray*}
\frac{1}{n\tau_{n}^{p}}\sum_{i=1}^{[nt]}\bigl(|\Delta
_{i}^{n}Y|^{p}-\bigl|\sigma
_{{(i-1)/n}}\Delta_{i}^{n}G\bigr|^{p}\bigr)\stackrel
{P}{\longrightarrow}0.
\end{eqnarray*}
By applying the inequality $||x|^{p}-|y|^{p}|\leq
p|x-y|(|x|^{p-1}+|y|^{p-1}) $ for $p>1$ and $||x|^{p}-|y|^{p}|\leq|x-y|^{p}$
for $p\leq1$, (\ref{incrin}) and the Cauchy--Schwarz inequality, we can
conclude that the above convergence follows from
%
\begin{equation}\label{conv1}
\frac{1}{n\tau_{n}^{2}}\sum_{i=1}^{[nt]}E\bigl[\bigl|\Delta_{i}^{n}Y-\sigma
_{{%
(i-1)/n}}\Delta_{i}^{n}G\bigr|^{2}\bigr]\rightarrow0.
\end{equation}
Observe the decomposition
\begin{eqnarray*}
\Delta_{i}^{n}Y-\sigma_{{(i-1)/n}}\Delta
_{i}^{n}G=A_{i}^{n}+B_{i}^{n,\varepsilon}+C_{i}^{n,\varepsilon} ,
\end{eqnarray*}
where
\begin{eqnarray*}
A_{i}^{n} &=&\int_{{(i-1)/n}}^{{i/n}}g\biggl({\frac
{i}{n}}-s\biggr)%
\bigl(\sigma_{s}-\sigma_{{(i-1)/n}}\bigr)W(\mathrm{d}s),
\\
B_{i}^{n,\varepsilon} &=&\int_{{{(i-1)/n}}-\varepsilon}^{
{(i-1)/n}}%
\biggl(g\biggl({\frac{i}{n}}-s\biggr)-g\biggl({\frac{i-1}{n}}-s\biggr)
\biggr)\sigma
_{s}W(\mathrm{d}s)
\\
&&{}-\sigma_{{(i-1)/n}}\int_{{
{(i-1)/n}}-\varepsilon
}^{{(i-1)/n}}\biggl(g\biggl({\frac{i}{n}}-s\biggr)-g\biggl({\frac
{i-1}{n}}-s\biggr)%
\biggr)W(\mathrm{d}s),
\\
C_{i}^{n,\varepsilon} &=&\int_{-\infty}^{{
{(i-1)/n}}-\varepsilon}\biggl(%
g\biggl({\frac{i}{n}}-s\biggr)-g\biggl({\frac{i-1}{n}}-s\biggr)
\biggr)\sigma_{s}W(%
\mathrm{d}s)
\\
&&{}-\sigma_{{(i-1)/n}}\int_{-\infty}^{{{(i-1)/n}}
-\varepsilon}\biggl(g\biggl({\frac{i}{n}}-s\biggr)-g\biggl({\frac
{i-1}{n}}-s\biggr)%
\biggr)W(\mathrm{d}s).
\end{eqnarray*}
By Lemma \ref{lemma2} and the boundedness of $\sigma$ on compact intervals,
we deduce
%
\begin{equation}\label{ineq1}
\frac{1}{n\tau_{n}^{2}}\sum_{i=1}^{[nt]}E[|C_{i}^{n,\varepsilon
}|^{2}]\leq
C\pi^{n}((\varepsilon,\infty)),
\end{equation}
and by (\ref{measurecond}) we obtain that
\begin{eqnarray*}
\lim_{n\rightarrow\infty}\frac{1}{n\tau_{n}^{2}}%
\sum_{i=1}^{[nt]}E[|C_{i}^{n,\varepsilon}|^{2}]=0.
\end{eqnarray*}
Next, we get
%
\begin{equation}\label{ineq2}
\frac{1}{n\tau_{n}^{2}}\sum_{i=1}^{[nt]}E[|A_{i}^{n}|^{2}]\leq\frac
{C}{%
n\tau_{n}^{2}}E\Biggl[\sum_{i=1}^{[nt]}\int_{{(i-1)/n}}^{
{i/n}}g^{2}%
\biggl({\frac{i}{n}}-s\biggr)\bigl(\sigma_{s}-\sigma_{
{(i-1)/n}}\bigr)^{2}\,\mathrm{d}s%
\Biggr].
\end{equation}
Set $v(s,\eta)=\sup\{|\sigma_{s}-\sigma_{r}|^{2}| s,r\in\lbrack
-t,t], |r-s|\leq\eta\}$. Then we obtain
%
\begin{equation}\label{ineq3}
\frac{1}{n\tau_{n}^{2}}\sum_{i=1}^{[nt]}E[|A_{i}^{n}|^{2}]\leq\frac
{1}{n}%
\sum_{i=1}^{[nt]}E\biggl[v\biggl(\frac{i-1}{n},n^{-1}\biggr)\biggr].
\end{equation}
Moreover, for any $\kappa>0,$ since $\sigma$ is cadlag, there
exists $n$ big enough\ that
\begin{eqnarray*}
v\biggl(\frac{i-1}{n},n^{-1}\biggr)\leq\kappa+\bigl( \Delta\sigma_{
{(i-1)/n}%
}\bigr) ^{2}\mathbf{1}_{\{ ( \Delta\sigma_{
{(i-1)/n}%
}) ^{2}\geq\kappa\} },
\end{eqnarray*}
so%
\begin{eqnarray*}
\frac{1}{n\tau_{n}^{2}}\sum_{i=1}^{[nt]}E[|A_{i}^{n}|^{2}] &\leq
&\kappa+%
\frac{1}{n}\sum_{i=1}^{[nt]}E\bigl[\bigl( \Delta\sigma_{
{(i-1)/n}%
}\bigr) ^{2}\mathbf{1}_{\{ ( \Delta\sigma_{
{(i-1)/n}%
}) ^{2}\geq\kappa\} }\bigr]
 \\
&\leq&\kappa+E\biggl[\frac{1}{n}\sum_{-t\leq s\leq t}^{{}}(
\Delta
\sigma_{s}) ^{2}\mathbf{1}_{\{ ( \Delta\sigma
_{s})
^{2}\geq\kappa\} }\biggr],
\end{eqnarray*}
then%
\begin{eqnarray*}
\lim_{n\rightarrow\infty}\frac{1}{n\tau_{n}^{2}}%
\sum_{i=1}^{[nt]}E[|A_{i}^{n}|^{2}]\leq\kappa
\end{eqnarray*}
and the convergence to zero follows, letting $\kappa$ tend to zero.

Finally, observe the decomposition $B_{i}^{n,\varepsilon
}=B_{i}^{n,\varepsilon}(1)+B_{i}^{n,\varepsilon}(2)$ with
\begin{eqnarray*}
B_{i}^{n,\varepsilon}(1) &=&\int_{{{(i-1)/n}}-\varepsilon
}^{{(i-1)/
n}}\biggl(g\biggl({\frac{i}{n}}-s\biggr)-g\biggl({\frac{i-1}{n}}-s
\biggr)\biggr)\bigl(\sigma
_{s}-\sigma_{{{(i-1)/n}}-\varepsilon}\bigr)W(\mathrm{d}s),
\\
B_{i}^{n,\varepsilon}(2) &=&\bigl(\sigma_{{{(i-1)/n}}-\varepsilon
}-\sigma
_{{(i-1)/n}}\bigr)\int_{{{(i-1)/n}}-\varepsilon}^{
{(i-1)/n}}\biggl(g%
\biggl({\frac{i}{n}}-s\biggr)-g\biggl({\frac{i-1}{n}}-s\biggr)
\biggr)W(\mathrm{d}s).
\end{eqnarray*}
We obtain the inequalities
\begin{eqnarray}\label{ineq4}
\frac{1}{n\tau_{n}^{2}}\sum_{i=1}^{[nt]}E[|B_{i}^{n,\varepsilon
}(1)|^{2}]&\leq&\frac{1}{n}\sum_{i=1}^{[nt]}E\biggl[v\biggl(\frac
{i-1}{n},\varepsilon
\biggr)\biggr],\nonumber
\\[-8pt]\\[-8pt]
\frac{1}{n\tau_{n}^{2}}\sum_{i=1}^{[nt]}E[|B_{i}^{n,\varepsilon
}(2)|^{2}]&\leq&\frac{1}{n}\sum_{i=1}^{[nt]}E\biggl[v\biggl(\frac
{i-1}{n},\varepsilon
\biggr)^{2}\biggr]^{1/2}. \nonumber
\end{eqnarray}
By using the same arguments as above, we have that both terms converge to
zero and we obtain (\ref{conv1}), which completes the proof with $Z_{1}=0.$

To prove the general case, \ with $Z_{1}\neq0$, we consider, for
simplicity, the case $k=2$.

Assume first that $0\leq p_{1},p_{2}\leq1$. We have $%
||x_{1}+y_{1}|^{p_{1}}|x_{2}+y_{2}|^{p_{2}}-|y_{1}|^{p_{1}}|y_{2}|^{p_{2}}|
\leq
C(|x_{1}|^{p_{1}}|x_{2}|^{p_{2}}+|x_{1}|^{p_{1}}|y_{2}|^{p_{2}}+|y_{1}|^{p_{1}}|x_{2}|^{p_{2}})
$. Hence we deduce
\begin{eqnarray*}
&&|V(Z,p_{1},p_{2})_{t}^{n}-V(Y,p_{1},p_{2})_{t}^{n}|
\\
&&\quad \leq\frac{C}{n\tau_{n}^{p_{+}}}\sum_{i=1}^{[nt]-1} ( |\Delta_i^n
Z_1|^{p_1} |\Delta_{i+1}^n Z_1|^{p_2} +|\Delta_i^n Z_1|^{p_1}
|\Delta_{i+1}^n Z_2|^{p_2}+|\Delta_{i}^n Z_2|^{p_1} |\Delta_{i+1}^n
Z_1|^{p_2} ) ,
\end{eqnarray*}
and the result follows by (\ref{pr1cond}).

Next, assume that $p_{1}\leq p_{2}$, $p_{2}>1$. We deduce that
\begin{eqnarray*}
&&|(V(Z,p_{1},p_{2})_{t}^{n})^{1/p_{2}}-(V(Y,p_{1},p_{2})_{t}^{n})^{1/p_{2}}|
\\
&&\quad\leq C\Biggl(\Biggl(\frac{1}{n\tau_{n}^{p_{+}}}\sum_{i=1}^{[nt]-1}
|\Delta_i^n Z_1|^{p_1} |\Delta_{i+1}^n Z_1|^{p_2}\Biggr)^{1/p_{2}}+
\Biggl(\frac{%
1}{n\tau_{n}^{p_{+}}}\sum_{i=1}^{[nt]-1} |\Delta_i^n Z_1|^{p_1}
|\Delta_{i+1}^n Z_2|^{p_2}\Biggr)^{1/p_{2}}
\\
&&{}\hspace*{13pt}\qquad+ \Biggl(\frac{1}{n\tau_{n}^{p_{+}}}\sum_{i=1}^{[nt]-1} |\Delta_i^n
Z_2|^{p_1} |\Delta_{i+1}^n Z_1|^{p_2}\Biggr)^{1/p_{2}} \Biggr) ,
\end{eqnarray*}
which completes the proof of Theorem \ref{th3}.
\end{pf*}

\subsection{Proof of Theorems 2 and 5}\label{sec8.5}

\begin{pf*}{Proof of Theorem \ref{th2}} We first show the weak
convergence of
finite-dimensional distributions and then prove the tightness of the
sequence $\sqrt{n}(V(p_{1}^{j},\ldots,p_{k}^{j})_{t}^{n}-\rho
_{p_{1}^{j},\ldots,p_{k}^{j}}^{(n)}t)_{1\leq j\leq d}$.

\textit{Step} 1: Define the vector $Z_{n}(j)=(Z_{n}^{1}(j),\ldots
,Z_{n}^{\mathrm{e}}(j))^{T}$, $1\leq j\leq d$, by
%
\begin{equation}\label{pry}
Z_{n}^{l}(j)={\frac{1}{\sqrt{n}}}\sum_{i=[nc_{l}]+1}^{[nb_{l}]}\bigl(
|X_{i,n}|^{p_{1}^{j}}\cdots|X_{i+k-1,n}|^{p_{k}^{j}}-\rho
_{p_{1}^{j},\ldots,p_{k}^{j}}^{(n)}\bigr) ,
\end{equation}
where $(c_{l},b_{l}]$, $l=1,\ldots,\mathrm{e}$, are disjoint intervals
contained in $%
[0,T]$. Set $Z_{n}^{l}=(Z_{n}^{l}(1),\ldots,\break Z_{n}^{l}(d))$,
$l=1,\ldots,\mathrm{e}$%
. Clearly, it suffices to prove that
\begin{eqnarray*}
( Z_{n}^{l}) _{1\leq l\leq \mathrm{e}}\stackrel{\mathcal{D}}{%
\longrightarrow}\bigl(\beta^{1/2}(B_{b_{l}}-B_{c_{l}})\bigr)_{1\leq
l\leq\mathrm{e}} ,
\end{eqnarray*}
where the matrix $\beta$ is given in Theorem \ref{th2}. By (\ref{chaosrep}),
we have the representation
\begin{eqnarray*}
Z_{n}^{l}(j)=\sum_{m=2}^{\infty}I_{k}\Biggl(\frac{1}{\sqrt{n}}%
\sum_{i=[nc_{l}]+1}^{[nb_{l}]}f_{m,j}^{n}(i)\Biggr).
\end{eqnarray*}
Set $F_{m,l}^{n}(j)=\frac{1}{\sqrt{n}}%
\sum_{i=[nc_{l}]+1}^{[nb_{l}]}f_{m,j}^{n}(i)$. By Theorem 2 in \cite{BNCPW08}
we obtain the weak convergence of finite-dimensional distributions when we
show that:

\begin{enumerate}[(iii)]
\item[(i)] For any $1\leq l\leq \mathrm{e}$, $1\leq j\leq d$, we have
\begin{eqnarray*}
\lim_{N\rightarrow\infty} \limsup_{n\rightarrow\infty}
\sum_{m=N+1}^\infty m!\|F_{m,l}^n (j)\|^2_{\mathbb{H}^{\otimes k}} =0.
\end{eqnarray*}

\item[(ii)] For any $m\geq2$, $1\leq l\leq \mathrm{e}$ and $1\leq j_1,j_2\leq
d$, we
have constants $C_{k,l}$ such that
\begin{eqnarray*}
\lim_{n\rightarrow\infty} m!\langle F_{m,l}^n (j_1),F_{m,l}^n
(j_2)\rangle_{\mathbb{H}^{\otimes m}} = C_{m,l} (j_1,j_2) ,
\end{eqnarray*}
and $\sum_{m=2}^\infty C_{m,l} (j_1,j_2)= \beta_{j_1, j_2} (b_l-c_l) $.

\item[(iii)] For any $1\leq l_1 \not= l_2 \leq \mathrm{e}$ and $1\leq
j_1,j_2\leq d$,
we have
\begin{eqnarray*}
\lim_{n\rightarrow\infty} \langle F_{m,l_1}^n (j_1),F_{m,l_2}^n
(j_2)\rangle_{\mathbb{H}^{\otimes m}} = 0.
\end{eqnarray*}

\item[(iv)] For any $m\geq2$, $1\leq l\leq \mathrm{e}$, $1\leq j\leq d$ and $p=1,
\ldots, m-1$
\begin{eqnarray*}
\lim_{n\rightarrow\infty} \|F_{m,l}^n (j) \otimes_p F_{m,l}^n
(j)\|^2_{\mathbb{H}^{\otimes2(m-p)}} = 0.
\end{eqnarray*}
\end{enumerate}

Note that it is sufficient to prove (i), (ii) and (iv) for $l=1$, $b_{l}=1$
and $a_{l}=0$. In this case we use the notation
$F_{m}^{n}(j)=F_{m,1}^{n}(j)$.

(i) and (ii) As in (\ref{impest}), we have
%
\begin{eqnarray}\label{estimate}
&&m!\langle F_{m}^{n}(j_{1}),F_{m}^{n}(j_{2})\rangle_{\mathbb
{H}^{\otimes m}}\nonumber
\\
&&\quad=m!%
\Biggl(\langle f_{m,j_{1}}^{n}(1),f_{m,j_{2}}^{n}(1)\rangle_{\mathbb
{H}^{\otimes
m}}+\frac{2}{n}\sum_{l=1}^{n-1}(n-l)\langle
f_{m,j_{1}}^{n}(1),f_{m,j_{2}}^{n}(1+l)\rangle_{\mathbb{H}^{\otimes
m}}\Biggr)\qquad
\\
&&\quad\leq Cm!\langle f_{m,j_{1}}^{n}(1),f_{m,j_{2}}^{n}(1)\rangle
_{\mathbb{H}^{\otimes m}}\Biggl(1+\sum_{l=1}^{n-1}\mathrm{r}(l)\Biggr).\nonumber
\end{eqnarray}
Since $\sum_{l=1}^{\infty}r(l)<\infty$, we obtain by (\ref
{converg1})--(\ref%
{converg3}) and the dominated convergence theorem
\begin{eqnarray*}
&&\lim_{n\rightarrow\infty}m!\langle
F_{m,l}^{n}(j_{1}),F_{m,l}^{n}(j_{2})\rangle_{\mathbb{H}^{\otimes
m}}
\\
&&\quad =C_{m}(j_{1},j_{2})
\\
&&\quad= m!\Biggl(\langle f_{m,j_{1}}(1),f_{m,j_{2}}(1)\rangle_{\mathbb
{H}^{\otimes
m}}+2\sum_{l=1}^{\infty}\langle
f_{m,j_{1}}^{n}(1),f_{m,j_{2}}^{n}(1+l)\rangle_{\mathbb{H}^{\otimes
m}}\Biggr) ,
\end{eqnarray*}
and $\sum_{m=2}^{\infty}C_{m}(j_{1},j_{2})=\beta_{j_{1},j_{2}}$ (notice
that $\beta_{j_{1},j_{2}}$ is finite due to the dominated convergence
theorem). Hence, we deduce (ii). On the other hand, we have
\begin{eqnarray*}
\limsup_{n\rightarrow\infty}\sum_{m=N+1}^{\infty}m!\langle
f_{m,j_{1}}^{n}(1),f_{m,j_{2}}^{n}(1)\rangle_{\mathbb{H}^{\otimes
m}}=\sum_{m=N+1}^{\infty}m!\langle
f_{m,j_{1}}(1),f_{m,j_{2}}(1)\rangle
_{\mathbb{H}^{\otimes m}}<\infty.
\end{eqnarray*}
Thus, we obtain (i) by (\ref{estimate}).

(iii) W.l.o.g. consider the case $j=j_{1}=j_{2}$. For any
$l_{1}<l_{2}$, as
in (\ref{impest}), we have
\begin{eqnarray*}
|\langle F_{m,l_{1}}^{n}(j),F_{m,l_{2}}^{n}(j)\rangle_{\mathbb
{H}^{\otimes
m}}|\leq\frac{C}{n}\sum_{h=[nc_{l_{1}}\hspace*{-1pt}]+1}^{[nb_{l_{1}}\hspace*{-1pt}]}%
\sum_{i=[nc_{l_{2}}\hspace*{-1pt}]+1}^{[nb_{l_{2}}\hspace*{-1pt}]}|r_{n}^{m}(i-h)|.
\end{eqnarray*}
Assume w.l.o.g. that $c_{l_{1}}=0$, $b_{l_{1}}=c_{l_{2}}=1$ and $b_{l_{2}}=2$
(the case $b_{l_{1}}<c_{l_{2}}$ is much easier). Then, by condition
(\ref%
{cltcond2}), we obtain the approximation (as in (\ref{estimate}))
\begin{eqnarray*}
|\langle F_{m,l_{1}}^{n}(j),F_{m,l_{2}}^{n}(j)\rangle_{\mathbb
{H}^{\otimes
m}}|\leq C\Biggl(\frac{1}{n}\sum_{j=1}^{n}j\mathrm{r}(j)+\sum
_{j=1}^{n-1}%
\mathrm{r}(n+j)\Biggr)\rightarrow0 ,
\end{eqnarray*}
since $\sum_{j=1}^{\infty}\mathrm{r}(j)<\infty$.

(iv) A straightforward computation shows that
\begin{eqnarray*}
&&\|F_{m}^{n}(j)\otimes_{p}F_{m}^{n}(j)\|_{\mathbb{H}^{\otimes
2(m-p)}}^{2}
\\
&&\quad={%
\frac{C_{m,j}}{n^{2}}}%
\sum
_{i_{1},i_{2},i_{3},i_{4}=1}^{n-1}r_{n}^{p}(|i_{1}-i_{2}|)r_{n}^{p}(|i_{4}-i_{3}|)r_{n}^{m-p}(|i_{1}-i_{4}|)r_{n}^{m-p}(|i_{2}-j_{3}|)
\end{eqnarray*}
for some constant $C_{m,j}$. The latter is smaller than
\begin{eqnarray*}
&&{\frac{C}{n}}%
\sum
_{i,h,l=1}^{n-1}|r_{n}^{p}(i)||r_{n}^{p}(l)|\bigl|r_{n}^{m-p}(|i-h|)\bigr|\bigl|r_{n}^{m-p}(|l-h|)\bigr|
\\
&&\quad=
{\frac{C}{n}}\sum_{h=1}^{n-1}\Biggl(%
\sum_{i=1}^{n-1}|r_{n}^{p}(i)|\bigl|r_{n}^{m-p}(|i-h|)\bigr|\Biggr)^{2}.
\end{eqnarray*}
Now, for any $0<\varepsilon<1$, we obtain by the H\"{o}lder inequality
\begin{eqnarray*}
&&n^{-1}\sum_{0\leq h\leq n-1}\biggl( \sum_{0\leq i\leq
n-1}|r_{n}^{p}(i)|\bigl|r_{n}^{m-p}(|i-h|)\bigr|\biggr) ^{2}
\\
&&\quad\leq n^{-1}\sum
_{0\leq
h\leq\lbrack n\varepsilon]}\biggl( \sum_{0\leq i\leq
n-1}|r_{n}^{p}(i)|\bigl|r_{n}^{m-p}(|i-h|)\bigr|\biggr) ^{2}
\\
&&{}\qquad+2n^{-1}\sum_{h=[n\varepsilon]}^{n-1}\Biggl( \sum
_{i=0}^{[n\varepsilon
/2]}|r_{n}^{p}(i)|\bigl|r_{n}^{m-p}(|i-h|)\bigr|\Biggr)^{2}
\\
&&{}\qquad+2n^{-1}\sum_{h=[n\varepsilon]}^{n-1}\Biggl( \sum
_{h=[n\varepsilon
/2]}^{n-1}|r_{n}^{p}(i)|\bigl|r_{n}^{m-p}(|i-h|)\bigr|\Biggr) ^{2}
\\
&&\quad\leq C\biggl( \varepsilon\biggl( \sum_{0\leq i\leq
n-1}|r_{n}^{m}(i)|\biggr) ^{2}+\biggl(\sum_{0\leq i\leq
n-1}|r_{n}^{m}(i)|\biggr)%
^{2p/m}\biggl(\sum_{[n\varepsilon/2]<h\leq n-1}|r_{n}^{m}(h)|\biggr)%
^{2(m-p)/m}\biggr).
\end{eqnarray*}
The latter is smaller (again by (\ref{cltcond2})) than
\begin{eqnarray*}
C\biggl( \varepsilon\biggl( \sum_{0\leq i\leq n-1}\mathrm
{r}(i)\biggr) ^{2}+%
\biggl(\sum_{0\leq i\leq n-1}\mathrm{r}(i)\biggr)^{2p/m}\biggl(\sum
_{[n\varepsilon
/2]<h\leq n-1}\mathrm{r}(h)\biggr)^{2(m-p)/m}\biggr)
\end{eqnarray*}
that converges to $C\varepsilon( \sum_{i=0}^{\infty}r(i)
) ^{2}$
as $n\rightarrow\infty$. Thus, we obtain (iv) by letting $\varepsilon
\rightarrow0$.

\textit{Step} 2: Clearly, it suffices to consider the case $d=1$, $%
p_{l}^{1}=p_{l}$. Set
\begin{eqnarray*}
\sqrt{n}\bigl(V(p_{1},\ldots,p_{k})_{t}^{n}-\rho_{p_{1},\ldots
,p_{k}}^{(n)}t\bigr)&=&\sum_{m=2}^{\infty}I_{m}\Biggl(\frac{1}{\sqrt
{n}}%
\sum_{i=1}^{[nt]}f_{m}^{n}(i)\Biggr)+\mathrm{O}_{p}(n^{-1/2})
\\
&=:&Z_{t}^{n}+\mathrm{O}(n^{-1/2})
\end{eqnarray*}
(where the approximation holds locally uniformly in $t$) and
\begin{eqnarray*}
Z_{t}^{n,N}=\sum_{m=2}^{N}I_{m}\Biggl({\frac{1}{\sqrt{n}}}%
\sum_{i=1}^{[nt]}f_{m}^{n}(i)\Biggr).
\end{eqnarray*}
In step 1, we have proved that conditions (i)--(iii) of Theorem 2 in
\cite%
{BNCPW08} are satisfied. Then by (2.3) of Theorem 2 in \cite{BNCPW08} and
the Cauchy--Schwarz inequality, we obtain the approximation
\begin{eqnarray*}
&&P(|Z_{t}^{n,N}-Z_{t_{1}}^{n,N}|\geq\lambda
,|Z_{t_{2}}^{n,N}-Z_{t}^{n,N}|\geq\lambda)
\\
&&\quad\leq{\frac{%
E^{1/2}[|Z_{t}^{n,N}-Z_{t_{1}}^{n,N}|^{4}]E^{1/2}[|Z_{t_{2}}^{n,N}-Z_{t}^{n,N}|^{4}]
}{\lambda^{4}}}
\\
&&\quad\leq C{\frac{\beta^{2}([nt]-[nt_{1}])([nt_{2}]-[nt])}{n^{2}\lambda
^{4}}}%
\leq C{\frac{\beta^{2}(t_{2}-t_{1})^{2}}{\lambda^{4}}}
\end{eqnarray*}
for any $t_{1}\leq t\leq t_{2}$ and $\lambda>0$. On the other hand,
(\ref%
{converg3}) and (\ref{estimate}) imply that
\begin{eqnarray*}
\lim_{N\rightarrow\infty}E[|Z_{t}^{n}-Z_{t}^{n,N}|^{2}]=0
\end{eqnarray*}
for any $n$ and any $t$. Using this we conclude that
\begin{eqnarray*}
P(|Z_{t}^{n}-Z_{t_{1}}^{n}|\geq\lambda
,|Z_{t_{2}}^{n}-Z_{t}^{n}|\geq
\lambda)\leq C{\frac{\beta^{2}(t_{2}-t_{1})^{2}}{\lambda^{4}}}
\end{eqnarray*}
for any $t_{1}\leq t\leq t_{2}$ and $\lambda>0$, from which we deduce the
tightness of the sequence $Z_{t}^{n}$ by Theorem 15.6 in \cite{Bil68}. This
completes the proof of Theorem \ref{th2}.
\end{pf*}

\begin{pf*}{Proof of Theorem \ref{th4}} We only consider the case
$Z_{1}=0$ (if
$Z_{1}\neq0$, we can proceed as in Theorem \ref{th3}). Also, for the
sake of simplicity, we take $d=1$, $k=1$ and $p_{1}=p$. We use the
decomposition from the proof of Theorem \ref{th3}:
%
\begin{eqnarray}\label{decclt}
&&\sqrt{n}\biggl(V(Y,p)_{t}^{n}-\mu_{p}\int_{0}^{t}|\sigma
_{s}|^{p}\,\mathrm{d}s\biggr)\nonumber
\\
&&\quad=
\sqrt{n}\Biggl({\frac{1}{n\tau_{n}^{p}}}\sum_{j=1}^{[lt]}\bigl|\sigma
_{{(j-1)/
l}}\bigr|^{p}\sum_{i\in I_{l}(j)}|\Delta_{i}^{n}G|^{p}-\mu
_{p}l^{-1}\sum_{j=1}^{[lt]}\bigl|\sigma_{{(j-1)/l}}\bigr|^{p}\Biggr)
 \\
&&{}\qquad+\frac{1}{\sqrt{n}\tau_{n}^{p}}\sum_{i=1}^{[nt]}\bigl(|\Delta
_{i}^{n}Y|^{p}-\bigl|\sigma_{{(i-1)/n}}\Delta_{i}^{n}G\bigr|^{p}
\bigr)+\overline{R%
}_{t}^{n,l} \nonumber
\end{eqnarray}
for any $l\leq n$, with
\begin{eqnarray*}
\overline{R}_{t}^{n,l} &=&{\frac{1}{\sqrt{n}\tau_{n}^{p}}}\Biggl(%
\sum_{i=1}^{[nt]}\bigl|\sigma_{{(i-1)/n}}\Delta
_{i}^{n}G\bigr|^{p}-\sum_{j=1}^{[lt]}\bigl|\sigma_{{(j-1)/l}}\bigr|^{p}\sum
_{i\in
I_{l}(j)}|\Delta_{i}^{n}G|^{p}\Biggr)
\\
&&{}+\sqrt{n}\mu_{p}\Biggl(l^{-1}\sum_{j=1}^{[lt]}\bigl|\sigma_{
{(j-1)/l}%
}\bigr|^{p}-\int_{0}^{t}|\sigma_{s}|^{p}\, \mathrm{d}s\Biggr).
\end{eqnarray*}
Observe that under the assumption \textup{(CLT)} we obtain the weak
convergence
\begin{eqnarray*}
\sqrt{n}\bigl(V(G,p)_{t}^{n}-\mu_{p}t\bigr)\stackrel{\mathcal
{L}}{\rightarrow}%
\sqrt{\beta}B_{t}
\end{eqnarray*}
(see Theorem \ref{th2}). Since $E[G_{t}(V(G,p)_{t}^{n}-\mu_{p}t)]=0$ for
any $t>0$, because $G$ has a symmetric distribution, we deduce (by
Theorem 5
in \cite{BNCP07}) that
\begin{eqnarray*}
\bigl(G_{t},\sqrt{n}\bigl(V(G,p)_{t}^{n}-\mu_{p}t\bigr)\bigr)\stackrel
{\mathcal{L%
}}{\rightarrow}\bigl(G_{t},\sqrt{\beta}B_{t}\bigr).
\end{eqnarray*}
Now, an application of the condition $D^{\prime\prime}$ from
Proposition 2
in \cite{AE78} shows that
\begin{eqnarray*}
\sqrt{n}\bigl(V(G,p)_{t}^{n}-\mu_{p}t\bigr)\stackrel
{\mathrm{st}}{\longrightarrow}%
\sqrt{\beta}B_{t}.
\end{eqnarray*}
By the properties of stable convergence, it follows immediately that
\begin{eqnarray*}
\sqrt{n}\Biggl({\frac{1}{n\tau_{n}^{p}}}\sum_{j=1}^{[lt]}\bigl|\sigma
_{{(j-1)/
l}}\bigr|^{p}\sum_{i\in I_{l}(j)}|\Delta_{i}^{n}G|^{p}-\mu
_{p}l^{-1}\sum_{j=1}^{[lt]}\bigl|\sigma_{{(j-1)/l}}\bigr|^{p}
\Biggr)\stackrel{\mathrm{st}}{%
\longrightarrow}\sqrt{\beta}\sum_{j=1}^{[lt]}\bigl|\sigma_{
{(j-1)/l}%
}\bigr|^{p}\Delta_{j}^{l}B
\end{eqnarray*}
for any fixed $l$. On the other hand, we have
\begin{eqnarray*}
\sqrt{\beta}\sum_{j=1}^{[lt]}\bigl|\sigma_{{(j-1)/l}}\bigr|^{p}\Delta
_{j}^{l}B%
\stackrel{P}{\longrightarrow}\sqrt{\beta}\int_{0}^{t}|\sigma
_{s}|^{p}%
\,\mathrm{d}B_{s}
\end{eqnarray*}
as $l\rightarrow\infty$.

Now we need to prove that the other summands in the decomposition (\ref
{decclt}) are negligible. The negligibility of the term $\overline{R}%
_{t}^{n,l}$ is shown as in the proof of Theorem 7 in \cite{BNCP07} but by
using condition (\ref{smooth}) instead of the H\"{o}lder continuity of
index $%
\gamma$. So we are left to prove that
%
\begin{equation} \label{convclt}
\frac{1}{\sqrt{n}\tau_{n}^{p}}\sum_{i=1}^{[nt]}\bigl(|\Delta
_{i}^{n}Y|^{p}-\bigl|\sigma_{{(i-1)/n}}\Delta_{i}^{n}G\bigr|^{p}
\bigr)\stackrel{P}{%
\longrightarrow}0.
\end{equation}
By applying, for $p\geq1$, the inequality $||x|^{p}-|y|^{p}|\leq
p|x-y|(|x|^{p-1}+|y|^{p-1}),$ (\ref{incrin}) and the Cauchy--Schwarz
inequality, and, for $p\leq1,$ $||x|^{p}-|y|^{p}|\leq|x-y|^{p}$ and the
Jensen inequality, we have\vspace*{2pt}
\begin{eqnarray*}
\frac{1}{\sqrt{n}\tau_{n}^{p}}\sum_{i=1}^{[nt]}E\bigl||\Delta
_{i}^{n}Y|^{p}-\bigl|\sigma_{{(i-1)/n}}\Delta_{i}^{n}G\bigr|^{p}\bigr|\leq
\frac{%
1}{\sqrt{n}\tau_{n}^{p\wedge1}}\sum_{i=1}^{[nt]}\bigl( E\bigl|\Delta
_{i}^{n}Y-\sigma_{{(i-1)/n}}\Delta_{i}^{n}G\bigr|^{2}\bigr) ^{
{%
(p\wedge1)/2}}.\vspace*{2pt}
\end{eqnarray*}
Now we use a similar decomposition as presented in the proof of Theorem
\ref%
{th3}:\vspace*{2pt}
\begin{eqnarray*}
\Delta_{i}^{n}Y-\sigma_{{(i-1)/n}}\Delta
_{i}^{n}G=A_{i}^{n}+B_{i}^{n,\varepsilon
_{n}^{(1)}}+\sum_{j=1}^{l}C_{i}^{n,\varepsilon_{n}^{(j)},\varepsilon
_{n}^{(j+1)}},\vspace*{2pt}
\end{eqnarray*}
where $A_{i}^{n}$, $B_{i}^{n,\varepsilon_{n}^{(1)}}$ are defined as
above, $%
0<\varepsilon_{n}^{(1)}<\cdots<\varepsilon_{n}^{(l)}<\varepsilon
_{n}^{(l+1)}=\infty$ and\vspace*{2pt}
\begin{eqnarray*}
C_{i}^{n,\varepsilon_{n}^{(j)},\varepsilon_{n}^{(j+1)}} &=&\int
_{{{(i-1)%
/n}}-\varepsilon_{n}^{(j+1)}}^{{{(i-1)/n}}-\varepsilon
_{n}^{(j)}}\biggl(%
g\biggl({\frac{i}{n}}-s\biggr)-g\biggl({\frac{i-1}{n}}-s\biggr)
\biggr)\sigma_{s}W(%
\mathrm{d}s)
\\[2pt]
&&{}-\sigma_{{(i-1)/n}}\int_{{{(i-1)/n}}-\varepsilon
_{n}^{(j+1)}}^{{%
{(i-1)/n}}-\varepsilon_{n}^{(j)}}\biggl(g\biggl({\frac
{i}{n}}-s\biggr)-g\biggl({%
\frac{i-1}{n}}-s\biggr)\biggr)W(\mathrm{d}s).\vspace*{2pt}
\end{eqnarray*}
By the assumption \textup{(CLT)} and Lemma \ref{lemma2}, we obtain
the following
inequalities (since $\sigma$ is bounded on compact intervals)
\begin{eqnarray*}
\frac{1}{\sqrt{n}\tau_{n}^{p\wedge1}} \sum_{i=1}^{[nt]}(
E|A_{i}^{n}|^{2}) ^{{(p\wedge1)/2}}&\leq& Cn^{-\gamma(
p\wedge1) +{1/2}},
\\
\frac{1}{\sqrt{n}\tau_{n}^{p\wedge1}} \sum_{i=1}^{[nt]}\bigl(
E\bigl|B_{i}^{n,\varepsilon_{n}^{(1)}}\bigr|^{2}\bigr) ^{{(p\wedge
1)/2}}&\leq&
Cn^{1/2}\bigl|\varepsilon_{n}^{(1)}\bigr|^{\gamma( p\wedge1) },
\\
\frac{1}{\sqrt{n}\tau_{n}^{p\wedge1}} \sum_{i=1}^{[nt]}\bigl(
E\bigl|C_{i}^{n,\varepsilon_{n}^{(j)},\varepsilon_{n}^{(j+1)}}\bigr|^{2}
\bigr) ^{%
{(p\wedge1)/2}}&\leq& Cn^{1/2}\bigl|\varepsilon_{n}^{(j+1)}\bigr|^{\gamma
\bigl(
p\wedge1\bigr) }\big|\pi^{n}\bigl(\bigl[\varepsilon_{n}^{(j)},\varepsilon
_{n}^{(j+1)}\bigr)\bigr|^{{(p\wedge1)/2}} ,
\\
&&\hspace*{-155pt} j=1,\ldots,l-1,
\\
\frac{1}{\sqrt{n}\tau_{n}^{p\wedge1}} \sum_{i=1}^{[nt]}\bigl(
E\bigl|C_{i}^{n,\varepsilon_{n}^{(l)},\varepsilon_{n}^{(l+1)}}\bigr|^{2}
\bigr) ^{%
{(p\wedge1)/2}}&\leq& Cn^{1/2}\pi^{n}\bigl(\bigl(\varepsilon
_{n}^{(l)},\infty\bigr)\bigr)^{%
{(p\wedge1)/2}}.
\end{eqnarray*}

Then we deduce by \textup{(CLT)} and Lemma \ref{lemma3}
\begin{eqnarray*}
\frac{1}{\sqrt{n}\tau_{n}^{p}}\sum_{i=1}^{[nt]}\bigl|\Delta
_{i}^{n}Y-\sigma_{%
{(i-1)/n}}\Delta_{i}^{n}G\bigr|^{p}\stackrel{P}{\longrightarrow}0.
\end{eqnarray*}
This completes the proof of Theorem \ref{th4}.
\end{pf*}

\section*{Acknowledgements}
Ole E. Barndorff-Nielsen and Mark Podolskij acknowledge financial support
from CREATES, funded by the Danish National Research Foundation, and from
the Thiele Centre. The work of Jos\'{e} Manuel Corcuera is supported by the
MEC Grant No. MTM2009-08218.

\printhistory


\begin{thebibliography}{99}
\bibitem{AE78} Aldous, D.J. and Eagleson, G.K. (1978). On mixing and
stability of limit theorems. \textit{Ann.  Probab.} \textbf{6} 325--331.
\MR{0517416}

\bibitem{BNBas10} Barndorff-Nielsen, O.E. and Basse-O'Connor, A.
(2010). Quasi Ornstein--Uhlenbeck processes. \textit{Bernoulli}. To appear.

\bibitem{BNCP07} Barndorff-Nielsen, O.E., Corcuera, J.M. and
Podolskij, M. (2009). Power variation for Gaussian processes with stationary
increments. \textit{Stochastic Process. Appl.} \textbf{119} 1845--1865.
\MR{2519347}

\bibitem{BNCP09b} Barndorff-Nielsen, O.E., Corcuera, J.M. and
Podolskij, M. (2009). Multipower variation for Brownian semistationary
processes (full version). CREATES research paper 2009-21, Aarhus
Univ. Available at
\url{http://www.econ.au.dk/research/research-centres/creates/research/research-papers/research-papers-2009/}.

\bibitem{BNCPW08} Barndorff-Nielsen, O.E., Corcuera, J.M.,
Podolskij, M. and Woerner, J.H.C. (2009). Bipower variation for Gaussian
processes with stationary increments. \textit{J. Appl. Probab.} \textbf{46}
132--150.
\MR{2508510}

\bibitem{BNGJPS06} Barndorff-Nielsen, O.E., Graversen, S.E.,
Jacod, J., Podolskij, M. and Shephard, N. (2006). A central limit theorem
for realised power and bipower variations of continuous
semimartingales. In
\textit{From Stochastic
Calculus to Mathematical Finance. Festschrift in Honour of A.N. Shiryaev}
(Y. Kabanov, R.~Liptser and J. Stoyanov, eds.) 33--68. Heidelberg: Springer.
\MR{2233534}


\bibitem{BNSch04} Barndorff-Nielsen, O.E. and Schmiegel, J. (2004).
L\'{e}vy-based tempo-spatial modelling: With applications to turbulence.
\textit{Uspekhi Mat. NAUK} \textbf{59} 65--91.
\MR{2068843}

\bibitem{BNSch07a} Barndorff-Nielsen, O.E. and Schmiegel, J.
(2007). Ambit processes: With applications to turbulence and cancer growth.
In  \textit{Stochastic Analysis and Applications: The Abel
Symposium 2005} (F.E. Benth, G.D. Nunno, T. Linstr\o m, B. \O ksendal  and T. Zhang, eds.)
 93--124. Heidelberg: Springer.
\MR{2397785}

\bibitem{BNSch07b} Barndorff-Nielsen, O.E. and Schmiegel, J.
(2008a). A stochastic differential equation framework for the timewise
dynamics of turbulent velocities. \textit{Theory Probab.
Appl.} \textbf{52} 372--388.

\bibitem{BNSch08a} Barndorff-Nielsen, O.E. and Schmiegel, J.
(2008b): Time change, volatility and turbulence. In  \textit{Proceedings of the
Workshop on Mathematical Control Theory and Finance} (A. Sarychev, A.
Shiryaev, M. Guerra and M.D.R. Grossinho, eds.)  29--53.  Berlin:
Springer.
\MR{2484103}

\bibitem{BNSch08b} Barndorff-Nielsen, O.E. and Schmiegel, J.
(2008c). Time change and universality in turbulence. Research Report 2007-8.
Thiele Centre for Applied Mathematics in Natural Science. Unpublished
manuscript.

\bibitem{BNSch09} Barndorff-Nielsen, O.E. and Schmiegel, J. (2009).
Brownian semistationary processes and volatility/intermittency. In
 \textit{Advanced
Financial Modelling. Radon Series Comp. Appl. Math.} \textbf{8}
 (H. Albrecher, W. Rungaldier and W. Schachermeyer, eds.) 1--26.
Berlin: W. de Gruyter.
\MR{2648456}



\bibitem{BNS04a} Barndorff-Nielsen, O.E. and Shephard, N. (2004).
Power and bipower variation with stochastic volatility and jumps (with
discussion). \textit{J. Fin. Econometrics} \textbf{2} 1--48.

\bibitem{BNS04b} Barndorff-Nielsen, O.E. and Shephard, N. (2004).
Econometric analysis of realised covariation: High frequency covariance,
regression and correlation in financial economics. \textit{Econometrica}
\textbf{72} 885--925.
\MR{2051439}

\bibitem{BNS06} Barndorff-Nielsen, O.E. and Shephard, N. (2006).
Impact of jumps on returns and realised variances: Econometric analysis of
time-deformed L\'{e}vy processes. \textit{J. Econometrics} \textbf{131}
217--252.
\MR{2276000}

\bibitem{BNS07} Barndorff-Nielsen, O.E. and Shephard, N. (2007).
Variation, jumps, market frictions and high frequency data in financial
econometrics. In  \textit{Advances in Economics and Econometrics. Theory and Applications.
Ninth World Congress} (R. Blundell, T. Persson and W.K. Newey, eds.) 328--372.
Cambridge Univ. Press.

\bibitem{BNSW06} Barndorff-Nielsen, O.E. and Shephard, N. and
Winkel, M. (2006). Limit theorems for multipower variation in the presence
of jumps. \textit{Stochastic Process. Appl.} \textbf{116} 796--806.
\MR{2218336}

\bibitem{Beg07a} B\'{e}gyn, A. (2007). Asymptotic expansions and
central limit theorem for quadratic variations of Gaussian processes.
\textit{Bernoulli} \textbf{13} 712--753.
\MR{2348748}

\bibitem{Beg07b} B\'{e}gyn, A. (2007). Functional limit theorems
for generalized quadratic variations of Gaussian processes. \textit{Stochastic
Process. Appl.} \textbf{117} 1848--1869.
\MR{2437732}

\bibitem{Bil68} Billingsley, P. (1968). \textit{Convergence of
Probability Measures}. New York: Wiley.
\MR{0233396}



\bibitem{CL67} Cram\'{e}r, H. and Leadbetter, M.R. (1967).\textit{\
Stationary and Related Stochastic Processes}. New York: Wiley.



\bibitem{GuyLe89} Guyon, X. and Leon, J. (1989). Convergence en loi
des H-variation d'un processus gaussien stationaire. \textit{Ann. Inst.
H.
Poincar\'{e} Probab. Statist.} \textbf{25} 265--282.
\MR{1023952}

\bibitem{HoSu87} Ho, H.C. and Sun, T.C.  (1987). A central limit
theorem for non-instantaneous filters of a stationary Gaussian process.
\textit{J. Multivariate Anal.} \textbf{22} 144--155.
\MR{0890889}


\bibitem{Jac08a} Jacod, J. (2008). Asymptotic properties of
realized power variations and related functionals of semimartingales.
\textit{Stochastic Process. Appl.} \textbf{118} 517--559.
\MR{2394762}

\bibitem{Jac08b} Jacod, J. (2008). Statistics and high frequency
data. Lecture notes. Department of Statistics, Universit\'{e} Paris VI.

\bibitem{Jac08c} Jacod, J. (2008). Asymptotic properties of
realized power variations and related functionals of semimartingales:
Multipower variations. Working paper.


\bibitem{KiPo07} Kinnebrock, S. and Podolskij, M. (2008). A note on
the central limit theorem for bipower variation of general functions.
\textit{Stochastic Process. Appl.} \textbf{118} 1056--1070.
\MR{2418258}

\bibitem{Nab52} Nabeya, S. (1952). Absolute moments in 3-dimensional
normal distribution. \textit{Ann. Inst. Statist. Math.} \textbf{4} 15--30.
\MR{0052072}


\bibitem{Nu06} Nualart, D. (2006). \textit{The Malliavin Calculus and
Related Topics}, 2nd ed. Berlin: Springer.
\MR{2200233}

\bibitem{NuPec05} Nualart, D. and Peccati, G. (2005). Central limit
theorems for sequences of multiple stochastic integrals. \textit{Ann. Probab.}
\textbf{33} 177--193.
\MR{2118863}

\bibitem{NuO-L08} Nualart, D. and Ortiz-Latorre, S. (2008). Central
limit theorems for multiple stochastic integrals and Malliavin
calculus.
\textit{Stochastic Process. Appl.} \textbf{118} 614--628.
\MR{2394845}

\bibitem{PecTu05} Peccati, G. and Tudor, C.A. (2005). Gaussian
limits for vector-valued multiple stochastic integrals. In  \textit{Seminaire de Probabilites XXXVIII}
  (M. Emery, M.
Ledoux and M. Yor, eds.). \textit{Lecture Notes in Mathematics} \textbf{1857} 247--262. Berlin: Springer.
\MR{2126978}

%
%
%

\bibitem{Woe06} Woerner, J.H.C. (2006). Power and multipower
variation: Inference for high frequency data. In \textit{Stochastic Finance} (A.N. Shiryaev, M. do Ros\'
{a}rio Grossihno, P. Oliviera and M. Esquivel, eds.) 343--364.
Heidelberg: Springer.
\MR{2230770}


\bibitem{Woe08} Woerner, J.H.C. (2008). Volatility estimates for
high frequency data: Market microstructure noise versus fractional Brownian
motion models. Unpublished manuscript.
\end{thebibliography}
\end{document}